\documentclass[10pt,twoside,fleqn,a4paper]{article}

\usepackage{amssymb}
\usepackage{epsfig}

\newtheorem{theorem}{Theorem}[section]

\newtheorem{definition}{Definition}[section]
\newtheorem{algorithm}{Algorithm}[section]

\title{
  Partitioning Sparse Graphs using the Second Eigenvector of
  their Graph Laplacian
}

\author{
  David De Wit
}

\date{June 1991}

\begin{document}

\maketitle

\begin{abstract}
  Partitioning a graph into three pieces, with two of them large
  and connected, and the third a small ``separator'' set, is
  useful for improving the performance of a number of
  combinatorial algorithms.  This is done using the second
  eigenvector of a matrix defined solely in terms of the
  incidence matrix, called the graph Laplacian.  For sparse
  graphs, the eigenvector can be efficiently computed using the
  Lanczos algorithm.  This graph partitioning algorithm is
  extended to provide a complete hierarchical subdivision of the
  graph.  The method has been implemented and numerical results
  obtained both for simple test problems and for several grid
  graphs.
\end{abstract}

\pagebreak

%%%%%%%%%%%%%%%%%%%%%%%%%%%%%%%%%%%%%%%%%%%%%%%%%%%%%%%%%%%%%%%%%%%%%%%%

\tableofcontents
\pagebreak

\listoffigures
\listoftables
\pagebreak

%%%%%%%%%%%%%%%%%%%%%%%%%%%%%%%%%%%%%%%%%%%%%%%%%%%%%%%%%%%%%%%%%%%%%%%%

\section{Introduction}
\label{sec:Introduction}

This report is a detailing of the research carried out by myself during
the four months of Semester 1 (March-June) 1991, at the Department of
Mathematics, University of Queensland, Australia.  It was carried out
under the supervision of Dr David E Stewart, and was accredited as a
\#30 project, subject classification code MN881.

%%%%%%%%%%%%%%%%%%%%%%%%%%%%%%%%%%%%%%%%%%%%%%%%%%%%%%%%%%%%%%%%%%%%%%%%

\subsection{Origins and Acknowledgements}
\label{sec:Origins and Acknowledgements}

This work was motivated by the publication of a
paper \cite{PothenSimonLiou:90}, titled ``Partitioning Sparse Matrices
with Eigenvectors of Graphs'', by Alex Pothen, Horst D. Simon and
Kang-Pu Liou in the {\it SIAM Journal on Matrix Analysis and
Applications}, September 1990, 11(3):430--452, and emulates some of the
implementation contained therein.  Their paper is mainly based on some
work appearing in an earlier paper \cite{Fiedler:75}, titled ``A
Property of Eigenvectors of Nonnegative Symmetric Matrices and its
Application to Graph Theory'', by Miroslav Fiedler in the {\it
Czechoslovak Mathematical Journal}, 1975, 25(100):619-633, which
provides a substantial amount of the theory quoted in
\S\ref{sec:Theory}.

I implemented these results in \textsc{C}, interfacing with a library of C
data-structures and functions called \textsc{meschach} \cite{Stewart:91},
written by my supervisor, Dr David Stewart.  He is also responsible for
selected pieces of code that I have used (see
Appendix \ref{app:Code}), as well as large amounts of time spent
educating me in \textsc{unix}, \textsc{C}, algorithms, numerical linear algebra and
scientific writing.  Proof reading was also done by my wife, Nerida
Iwasiuk.

%%%%%%%%%%%%%%%%%%%%%%%%%%%%%%%%%%%%%%%%%%%%%%%%%%%%%%%%%%%%%%%%%%%%%%%%

\subsection{Motivation and Applications}
\label{sec:Motivation and Applications}

The problem dealt with in this report is to partition the set of
vertices $ N = \left\lbrace 1, \dots, n \right\rbrace $ of an
undirected, unvaluated graph $ G = (N, E) $ into $3$ disjoint sets; a
minimally small separator set $ S $, and $3$ ``banks'' $ A $ and $ B $,
of size of order $ n / 2 $ each.

Immediate applications of this partitioning apply the
``divide-and-conquer''\\
paradigm to a range of graph-theoretic problems
such as the travelling sales representative problem.  Industrial
application problems, in particular those relating to the layout of
components in VLSI design, are discussed in \cite{Liu:89}.  Other
application problems are available in the Boeing-Harwell sparse matrix
test problem library \cite{DuffGrimesLewis:89}, and include large-scale
network-analysis problems such as power-grid distribution problems.

An immediate application in numerical linear algebra is in the
efficient solution of large sparse linear systems via factorisation and
parallel solution of subproblems.  Methods for doing this are outlined
below.

\pagebreak

Consider the solution of the $ n \times n $ (sparse) symmetric linear
system $ C \mathbf{x} = \mathbf{f} $.  If a permutation of the indices
of $ C $ can be made, such that in block form:

\begin{eqnarray*}
  P C P^{\top} \mathbf{x}
  =
  \left[
    \begin{array}{ccc}
      H_A & . & M_A^{\top} \\
      . & H_B & M_B^{\top} \\
      M_A & M_B & H_S
    \end{array}
  \right]
  \left[
    \begin{array}{c}
      \mathbf{x}_A \\ \mathbf{x}_B \\ \mathbf{x}_S
    \end{array}
  \right]
  =
  \left[
    \begin{array}{c}
      \mathbf{f}_A \\ \mathbf{f}_B \\ \mathbf{f}_S
    \end{array}
  \right],
\end{eqnarray*}
then the system can often be solved much faster.  Here $ H_A $,
$ H_B $ and $ H_S $ are symmetric, and $ M_A $, $ M_B $ are
usually not (they are usually not even square).  The system can be
solved using the $ L D L^{\top} $ (Cholesky-type) factorisation for
indefinite, symmetric matrices.  Firstly, be aware that simple Cholesky
factorisation $ L L^{\top} $ is only applicable to positive definite
matrices, and problems dealt with are \textbf{not} guaranteed to be
positive definite.  (Note that most of the experimental graphs used in
this report are 5-point grid graphs, which are always positive
semi-definite.) Whilst the graph Laplacian matrices mentioned here are
singular, the $ C $ matrices actually involved are not.  Writing
$ P C P^{\top} = L D L^{\top} $ in block form gives:
\begin{eqnarray*}
  \left[
    \begin{array}{ccc}
      H_A & . & M_A^{\top} \\
      . & H_B & M_B^{\top} \\
      M_A & M_B & H_S
    \end{array}
  \right]
  =
  \left[
    \begin{array}{ccc}
      L_A & . & . \\
      . & L_B & . \\
      N_A & N_B & L_S
    \end{array}
  \right]
  D
  \left[
    \begin{array}{ccc}
      L_A^{\top} & . & N_A^{\top} \\
      . & L_B^{\top} & N_B^{\top} \\
      . & . & L_S^{\top}
    \end{array}
  \right],
\end{eqnarray*}
where $ L_A $, $ L_B $ are lower triangular, $ L_S $ is
symmetric but not generally triangular, and $ N_A $ and $ N_B $ are
generally sparse $ | S | \times | A | $ and
$ | S | \times | B | $ matrices.  In a similar procedure to the
usual Cholesky procedure, the solution of the system is done in stages
(3, not 2, as there is also a $ D $).  The first stage in the
solution of $ L ( D L^{\top} \mathbf{x} )  = \mathbf{f} $ is the
solution for $ \mathbf{y} = D L^{\top} \mathbf{x} $ of
$ L \mathbf{y} = \mathbf{f} $ by:
\begin{eqnarray*}
  \left[
    \begin{array}{ccc}
      L_A & . & . \\
      . & L_B & . \\
      N_A & N_B & L_S
    \end{array}
  \right]
  \left[
    \begin{array}{c}
      \mathbf{y}_A \\ \mathbf{y}_B \\ \mathbf{y}_S
    \end{array}
  \right]
  =
  \left[
    \begin{array}{c}
      \mathbf{f}_A \\ \mathbf{f}_B \\ \mathbf{f}_S
    \end{array}
  \right].
\end{eqnarray*}
To solve this block system, it can be written as:
\begin{eqnarray*}
  L_A \mathbf{y}_A
  & = &
  \mathbf{f}_A \\
  L_B \mathbf{y}_B
  & = &
  \mathbf{f}_B \\
  L_S \mathbf{y}_S
  & = &
  \mathbf{f}_S - N_A \mathbf{y}_A - N_B \mathbf{y}_B.
\end{eqnarray*}
The third equation is a (dense) order $ | S | $ system.  If $ | S | $
is small, it can be solved at no great expense, using simple Gaussian
elimination.  This solution depends on the prior solution of the two
earlier systems for $ \mathbf{y}_A $ and $ \mathbf{y}_B $.  These are
sparse, linear, lower triangular systems of size roughly $ n/2 $, and
are cheap to solve using forward substitution.

The second stage is the (trivial) solution of
$ D ( L^{\top} \mathbf{x} ) = \mathbf{y} $, using
$ L^{\top} \mathbf{x} = \mathbf{z} $ and $ D \mathbf{z} = \mathbf{y} $;
which is $ \mathbf{z} = D^{-1} \mathbf{y} $, as $ D $ is a diagonal
matrix.  The third stage is similar to the first stage, and provides a
solution of $ L^{\top} \mathbf{x} = \mathbf{z} $, which is the solution
to the whole problem $ C \mathbf{x} = \mathbf{f} $:
\begin{eqnarray*}
  \left[
    \begin{array}{ccc}
      L_A^{\top} & . & N_A^{\top} \\
      . & L_B^{\top} & N_B^{\top} \\
      . & . & L_S^{\top}
    \end{array}
  \right]
  \left[
    \begin{array}{c}
      \mathbf{x}_A \\ \mathbf{x}_B \\ \mathbf{x}_S
    \end{array}
  \right]
  =
  \left[
    \begin{array}{c}
      \mathbf{z}_A \\ \mathbf{z}_B \\ \mathbf{z}_S
    \end{array}
  \right].
\end{eqnarray*}
Similarly to the first stage, we solve first for $ \mathbf{x}_S $ the
order $ | S | $ dense system
$ L_S^{\top} \mathbf{x}_S = \mathbf{z}_S $, then substitute this into:
\begin{eqnarray*}
  L_A^{\top} \mathbf{x}_A
  & = &
  \mathbf{z}_A - N_A^{\top} \mathbf{x}_S \\
  L_B^{\top} \mathbf{x}_B
  & = &
  \mathbf{z}_B - N_B^{\top} \mathbf{x}_S,
\end{eqnarray*}
to obtain $ \mathbf{x}_A $ and $ \mathbf{x}_B $ cheaply, by backward
substitution.

An alternative to the $ L D L^{\top} $ factorisation is the direct
solution of $ ( P C P^{\top} ) \mathbf{x} = \mathbf{f} $.  This
method involves recursively partitioning the large sparse matrix $ C $,
and is called ``Nested Dissection'' \cite{George:73}.  We write out
the equations in block form:

\begin{eqnarray*}
  H_A \mathbf{x}_A
  & = &
  \mathbf{f}_A - M_A^{\top} \mathbf{x}_S \\
  H_B \mathbf{x}_B
  & = &
  \mathbf{f}_B - M_B^{\top} \mathbf{x}_S \\
  H_S \mathbf{x}_S
  & = &
  \mathbf{f}_S - M_A \mathbf{x}_A - M_B \mathbf{x}_B.
\end{eqnarray*}

``Solving'' the first $2$ equations for $ \mathbf{x}_A $ and
$ \mathbf{x}_B $, then substituting the results into the third yields
an order $ | S | $ dense system, which gives $ \mathbf{x}_S $ as the
solution to:
\begin{eqnarray*}
  M_A H_A^{-1}
  \mathbf{f}_A - M_A^{\top} \mathbf{x}_S
  +
  M_B H_B^{-1}
    \mathbf{f}_B - M_B^{\top} \mathbf{x}_S
  +
  H_S \mathbf{x}_S
  =
  \mathbf{f}_S.
\end{eqnarray*}
Explicitly:
\begin{eqnarray*}
  \hspace{-11pt}
  \mathbf{x}_S
  =
  {
    (
      M_A H_A^{-1} M_A^{\top} +
      M_B H_B^{-1} M_B^{\top} +
      H_S
    )
  }^{-1}
  (
    \mathbf{f}_S -
    M_A H_A^{-1} \mathbf{f}_A -
    M_B H_B^{-1} \mathbf{f}_B
  ).
\end{eqnarray*}
Note that whilst the graph Laplacian is singular, the original $ C $
matrix is generally non-singular, so there is no problem writing this
explicitly.  After the solution for $ \mathbf{x}_S $, we have
$\mathbf{x}_A $ and $ \mathbf{x}_B $ as the solution to the following
order $ | A | $ and $ | B | $ sparse systems, explicitly written as:
\begin{eqnarray*}
  \mathbf{x}_A
  & = &
  H_A^{-1} ( \mathbf{f}_A  - M_A^{\top} \mathbf{x}_S )
  \\
  \mathbf{x}_B
  & = &
  H_B^{-1} ( \mathbf{f}_B  - M_B^{\top} \mathbf{x}_S ).
\end{eqnarray*}

If $ | A | \; \approx \; | B | \; \approx n/2 $, and
$ | S | $ is small, the computational expense is minimised.  This
situation is amenable to implementation on parallel processing
machines, since the $2$ linear systems are decoupled.  If these are too
large to solve simply, the decomposition of $ C $ can be repeated on
$ H_A $ and $ H_B $.  The factorisation process can be recursively
continued until units of a desired ``atomic'' size are obtained --
perhaps somewhere between $3$ and $10$.  A more thorough analysis could
be performed to decide the optimal choice, such that overall
computational expense is minimised.  This decomposition has been
performed for several example graphs in \S\ref{sec:Results and
Problems}.

In this report, the graphs considered are sparse (the average degree of
the vertices is low), and large (at current computational capabilities,
$ n $ is $ 10^3 $ to $ 10^5 $).  \cite{Liu:89,PothenSimonLiou:90}
illustrate graphs with $ n $ up to $32 400$ (for a $ 180 \times 180 $
5-point grid, average degree 4).  An example of a 5-point grid graph is
provided in \S\ref{sec:Ishmail}.  (9-point grid graphs are
similar, except that each vertex is connected to its nearest 8
neighbours, rather than its nearest 4.) Larger scale applications are
available in the Boeing-Harwell sparse matrix test problem
library \cite{DuffGrimesLewis:89}, which quotes examples of size up to
$44 609$.

%%%%%%%%%%%%%%%%%%%%%%%%%%%%%%%%%%%%%%%%%%%%%%%%%%%%%%%%%%%%%%%%%%%%%%%%

\subsection{Broad Outline of the Report}

The algorithm for partitioning graphs is dependent on analysis of some
of the properties of ``M-Matrices'', of which the graph Laplacian is an
instance.  \S\ref{sec:M-Matrices} describes these entities
and some of their properties.  The graph Laplacian is introduced, and
shown to be an M-Matrix in \S\ref{sec:Graph Laplacian}.

An amazing result (of Fiedler, 1975 \cite{Fiedler:75}), presented in
Algorithm \ref{alg:EdgeSeparator} on
page \pageref{alg:EdgeSeparator}, relates the discrete
mathematics of the interconnections of the graph to the continuous
mathematics of the eigendecomposition of its ``graph Laplacian''
(defined in \S\ref{sec:Graph Laplacian} in terms of the incidence
matrix of the graph). This result takes the components of the second
eigenvector of the graph Laplacian matrix, and uses them in a valuation
of the \textbf{vertices} of the graph.  An edge separator set of the
graph is then computed, such that the vertices are partitioned by this
set into $2$ connected banks of very similar size.  Using this edge
separator set, a combinatorial algorithm creates the $3$ vertex sets
described in \S\ref{sec:Motivation and Applications}. This algorithm is
fairly efficient, and has a number of immediate applications. Many
large tasks can be simplified by using this algorithm in a
``divide-and-conquer'' approach, which can greatly improve
computational efficiency.  The Lanczos algorithm (to determine
eigenvalues of large, sparse matrices) is implemented as a part of the
overall algorithm, and is discussed in \S\ref{sec:LanczosAlgorithm}.

This report firstly details experiments performed in synthesising the
results of \cite{PothenSimonLiou:90}, who used the above strategy in
computing separator sets of sparse graphs of order $ 10^3 $ to
$ 10^5 $ vertices.  Secondly, it demonstrates the \textbf{recursive}
decomposition of a graph down to atomic-unit sized subgraphs.  This
appears to be a new (but probably obvious) implementation.

%%%%%%%%%%%%%%%%%%%%%%%%%%%%%%%%%%%%%%%%%%%%%%%%%%%%%%%%%%%%%%%%%%%%%%%%

\subsection{Contents of the Report}

\S\ref{sec:Theory} is the largest part of this report, and
contains seven subsections.  It describes the theory behind the entire
graph partitioning process.  \S\ref{sec:Implementation}
describes the implementation in \textsc{C} of the algorithms in
\S\ref{sec:Theory}, and references the source code in
Appendix \ref{app:Code}.
\S\ref{sec:Results and Problems} mentions the success of the
code in dealing with a number of test examples.  Results are compared,
where possible to literature and touchstone cases (in particular,
partitioning 5-point grid graphs). \S\ref{sec:Results and Problems}
also describes problems encountered in programming and
implementation.  \S\ref{sec:Further Work} outlines the
directions in which the work could be continued, to make the code more
useful.  The decomposition of a small example graph (Moshe) is
presented in Appendix \ref{sec:Moshe}.

\pagebreak

%%%%%%%%%%%%%%%%%%%%%%%%%%%%%%%%%%%%%%%%%%%%%%%%%%%%%%%%%%%%%%%%%%%%%%%%

\section{Theory}
\label{sec:Theory}

This section has seven subsections:
\begin{enumerate}
\item
  M-Matrices, their origins, applications, definitions, and some
  interesting theorems.
\item
  Graph Laplacian, its definition, and properties.
\item
  Graph Partitioning, using the second eigenvalue of the graph
  Laplacian \cite{Fiedler:75}.
\item
  Spectral Partitioning Algorithm, developed from the graph
  partitioning result \cite{PothenSimonLiou:90}.
\item
  Minimum Cover -- algorithms involved in finding the minimum
  vertex cover of a bipartite graph.
\item
  Lanczos Algorithm \cite{GolubVanLoan:89} -- required for the
  implementation of the spectral partitioning algorithm on sparse
  matrices.
\item
  Recursive Decomposition -- recursive implementation of the
  partitioning process.
\end{enumerate}

In this report, $ G = (N, E) $ refers to a graph, where
$ N = \left\lbrace 1, \dots, n \right\rbrace $ is the set of vertices.
$ A $ always refers to a real, square matrix of size $ n $.  This is a
sufficient, but not necessary requirement for some of the definitions
and theorems in \S\ref{sec:M-Matrices}.  In later sections, specific
theorems that also depend on $ A $ being \textbf{symmetric} are
quoted.

%%%%%%%%%%%%%%%%%%%%%%%%%%%%%%%%%%%%%%%%%%%%%%%%%%%%%%%%%%%%%%%%%%%%%%%%

\subsection{M-Matrices}
\label{sec:M-Matrices}

M-Matrices are a subclass of real matrices that are closely related to
non-negative matrices \cite{Minc:88}, and have a number of interesting
properties (for instance, see Theorem \ref{thm:M-Matrix
Non-Negative Inverse}). They arise in several specific fields, e.g.:

\begin{enumerate}
\item
  The discretisation of boundary-value partial differential equations
  (both symmetric and non-symmetric) generates matrices which are the
  negative of M-Matrices. Several examples of the discretisation of
  Poisson's Equation (viz $ {\nabla}^2 u(x,y) = f(x,y) $) on rectangular
  5-point grids are explicitly demonstrated in
  \S\ref{sec:Results and Problems} of this report. See also
  \S\ref{sec:Graph Laplacian}.
\item
  Continuous-Time Markov Processes. The differential equations for
  probabilities are usually (singular) M-Matrices.
\item
  Economics.
\end{enumerate}

M-Matrices were introduced in 1937 in \cite{Ostrowski:37}, and their
properties have been investigated by a number of researchers.  Material
presented in this section is abstracted from several
sources \cite{BermanPlemmons:79,Fiedler:75,Minc:88,PooleBoullion:74}.
M-Matrices have a multiplicity of possible definitions
(e.g. \cite{BermanPlemmons:79} lists $50$ definitions for nonsingular
M-Matrices).  This allows great flexibility in proving results which
involve them.  The following definition of an M-Matrix is quoted
from \cite{Minc:88}.

\begin{definition}
  $ A $ is an \textbf{M-Matrix} if there is a non-negative matrix
  $ B \in \mathbb{R}^{n \times n} $ with maximal eigenvalue $ r $,
  and $ c \geqslant r $ such that $ A = c I_n - B $.
\end{definition}

The main diagonal entries of an M-Matrix are non-negative and all of
its other entries are non-positive.  Fiedler \cite{Fiedler:75} calls an
M-Matrix a matrix of class $ K_0 $, and a nonsingular M-Matrix a
matrix of class $ K $.

\begin{definition}
  $
    F_n
    =
    \left\lbrace A \; | \; a_{ij} \leqslant 0, \; i \neq j \right\rbrace
  $.
\end{definition}

That is, $ F_n $ is the set of real matrices whose off-diagonal
elements are non-positive. The following results are some necessary and
sufficient conditions for an element of $ F_n $ to be an M-Matrix,
and are provided mainly as an introduction to some of the properties of
M-Matrices. The first is a sufficient condition:

\begin{theorem}
  $ A \in F_n $ is an M-Matrix iff all its eigenvalues have a
  non-negative real part.
\end{theorem}

The second is a necessary condition:

\begin{theorem}
  $ A \in F_n $ is an M-Matrix iff every real eigenvalue of $ A $
  is non-negative.
\end{theorem}

\begin{definition}
  Given $ M $ as a non-empty subset of $ N $, define the
  \textbf{principal submatrix} of $ A $ corresponding to $ M $ as
  $ A(M) = \left\lbrace a_{ij} | i, j \in M \right\rbrace $.
\end{definition}

\begin{theorem}
  A principal submatrix of an M-Matrix is an M-Matrix.
\end{theorem}

\begin{theorem}
  $ A \in F_n $ is an M-Matrix iff all its principal minors are
  non-negative.
\end{theorem}

The next result is important, and is mentioned in \cite{Fiedler:75}.

\begin{theorem}
  $ A \in F_n $ is a nonsingular M-Matrix iff $ A^{-1} $ is
  non-negative.
  \label{thm:M-Matrix Non-Negative Inverse}
\end{theorem}

Another definition of $ M $-matrices is:

\begin{definition}
  $ A $ is an \textbf{M-Matrix} iff all its off-diagonal entries are
  $\leqslant 0 $, and all its principal minors are non-negative. It is
  a \textbf{nonsingular M-Matrix} if its principal minors are
  positive.
  \label{def:M-M}
\end{definition}

Theorems \ref{thm:ZeroEvalminus1} and \ref{thm:ZeroEval} are proven
in \cite{Fiedler:62} and \cite{Varga:62}.

\begin{theorem}
  A nonsingular M-Matrix has $ A^{-1} $ non-negative, and if it is
  irreducible, $ A^{-1} $ is positive (that is, none of its
  entries are zero).
  \label{thm:ZeroEvalminus1}
\end{theorem}

\begin{theorem}
  If $ A $ is an irreducible, singular M-Matrix, then:

  \begin{enumerate}
  \item
    $ 0 $ is a simple eigenvalue of $ A $.
  \item
    There is, up to a scale factor, a unique non-zero
    eigenvector $ \mathbf{u} \in \mathbb{R}^n $ with eigenvalue
    zero, and all the components of $\mathbf{u}$ are non-positive or
    non-negative.
  \end{enumerate}
  \label{thm:ZeroEval}
\end{theorem}

Theorem \ref{thm:ZeroEvalplus1} applies to \textbf{symmetric}
matrices, and follows from Definition \ref{def:M-M} and the
properties of positive definite (or semidefinite) matrices.

\begin{theorem}
  A symmetric matrix is an M-Matrix if all its off-diagonal entries are
  non-positive and all its eigenvalues non-negative. If its eigenvalues
  are positive, then it is a non-singular M-Matrix.
  \label{thm:ZeroEvalplus1}
\end{theorem}

%%%%%%%%%%%%%%%%%%%%%%%%%%%%%%%%%%%%%%%%%%%%%%%%%%%%%%%%%%%%%%%%%%%%%%%%

\subsection{Graph Laplacian}
\label{sec:Graph Laplacian}

This section introduces the graph Laplacian, and defines it in a number
of ways.

\begin{definition}
  Given a graph $ G = (N, E) $, define the \textbf{degree} of vertex
  $ i $ as $ d_i $, the number of edges in $ E $ with one end
  being vertex $ i $.
\end{definition}

\begin{definition}
  Given a graph $ G = (N, E) $, the
  \textbf{graph Laplacian} \cite{AndersonMorley:71,Fiedler:75} is the
  matrix $ L(G) $ (here $ G $ is the incidence matrix of the graph)
  of the quadratic form:
  \begin{eqnarray*}
    \mathbf{x}^{\top} L(G) \mathbf{x}
    =
    (L(G) \mathbf{x}, \mathbf{x} )
    =
    {\displaystyle \sum_{(i, j) \in E} {( x_i - x_j)}^2}.
  \end{eqnarray*}
  Thus $ L(G) = (l_{ij}) $, where
  $
    l_{ij} = l_{ji} =
    \left\{
    \begin{array}{cl}
      0  & i \neq j, \; (i, j) \notin E \\
      -1  & i \neq j, \; (i, j) \in E \\
      d_i  & i = j.
    \end{array}
    \right.
  $
\end{definition}

The structure of the graph Laplacian is easier to see in another way.
Define $ M \in \mathbb{R}^{n \times n} $ as the \textbf{adjacency}
matrix of $ G $, and
$ D = \mathrm{diag}( d_1, \dots, d_n ) \in \mathbb{R}^{n \times n} $.
The graph Laplacian is then $ L(G) = D - M $.  This
definition means that $ L(G) $ is immediately seen as being a matrix
with $ -1 $ replacing $ 1 $ in $ M $, and a diagonal whose $ i $th
element is the number of non-zero off-diagonal elements in row or
column $ i $ of $ M $.  This is illustrated by Figures \ref{fig:Idit}
and \ref{fig:IditGrLap}, and is the form used in implementation.

As the row and column sums of the graph Laplacian are zero,
$ \mathbf{e}  = {[ 1, \dots, 1 ]}^{\top} \in \mathbb{R}^{n} $ is an
eigenvector corresponding to the eigenvalue 0.  As it has a zero
eigenvalue, the graph Laplacian is singular.  As it satisfies the
requirements of Theorem \ref{thm:ZeroEvalplus1}, it is an
M-Matrix.  Hence, the graph Laplacian is a singular M-Matrix, and is
able to be used in some of the theorems in
\S\ref{sec:Graph Partitioning}.

Experimental graphs considered in this report include several 5-point
grid graphs associated with the discretisation of elliptic
Boundary-Value PDEs on rectangular regions.  If the rectangle is
reduced to an $ m \times n $ grid, the graph Laplacian is an
$ mn \times mn $ sparse matrix, of very definite structure.  Consider
\textbf{Poisson's Equation} in two dimensions, with boundary conditions
for $ u ( x, y ) $ on a rectangle:

\begin{eqnarray*}
  \nabla^2 u(x,y)
  =
  {\displaystyle \frac{\partial^2 u(x,y)}{\partial x^2}} +
  {\displaystyle \frac{\partial^2 u(x,y)}{\partial y^2}} =
  f(x,y).
\end{eqnarray*}
(It is \textbf{Laplace's Equation} if $ f ( x, y ) = 0 $.)
Discretisation on a regular $ m \times n $ grid, using the notation
$f_{ij} = f(x_i, y_j ) $, and $ u_{ij} = u ( x_i, y_j ) $ gives:
\begin{eqnarray*}
  u_{i+1, j} + u_{i-1, j} + u_{i, j-1} + u_{i, j+1} - 4 u_{ij}
  =
  h^2 f_{ij} \; \; \;
  \left\{
    \begin{array}{l}
      i = 2, \dots, m-1 \\
      j = 2, \dots, n-1.
    \end{array}
  \right.
\end{eqnarray*}

This generates an $ mn \times mn $ linear system
$C\mathbf{u}=\mathbf{f} $ on a 5-point grid graph $ G $. The graph
Laplacian $ L(G) $ is an $ m \times m $ block-tridiagonal matrix, where
each block is a symmetric $ n \times n $ matrix.

\begin{eqnarray*}
  L(G)
  =
  \left[
    \begin{array}{*{6}{c}}
      P      & I_n    & .      & .      & \dots   & .      \\
      I_n    & P      & I_n    & .      & \dots   & .      \\
      .      & I_n    & P      & I_n    & \dots   & .      \\
      \vdots & \vdots & \ddots & \ddots & \ddots  & \vdots \\
      .      & .      & .      & I_n    & P       & I_n    \\
      .      & .      & .      & .      & I_n     & P
    \end{array}
  \right]
\end{eqnarray*}
where:
\begin{eqnarray*}
  P
  =
  \left[
    \begin{array}{*{6}{c}}
      -4     & 1      & .      & .      & \dots  & .      \\
      1      & -4     & 1      & .      & \dots  & .      \\
      .      & 1      & -4     & 1      & \dots  & .      \\
       \vdots & \vdots & \ddots & \ddots & \ddots & \vdots \\
       .      & .      & .      & 1      & -4     & 1      \\
      .      & .      & .      & .      & 1      & -4
    \end{array}
  \right].
\end{eqnarray*}

The structure of $ C $ is closely related to this. A specific example
is presented in Appendix \ref{sec:Ishmail}. For 9-point grids or
non-elliptic PDEs, the structures are different, but the solution
techniques are essentially the same.

%%%%%%%%%%%%%%%%%%%%%%%%%%%%%%%%%%%%%%%%%%%%%%%%%%%%%%%%%%%%%%%%%%%%%%%%

\subsection{Graph Partitioning}
\label{sec:Graph Partitioning}

This section outlines the work of \cite{Fiedler:75}, and theorems
listed are taken directly from it.

\begin{definition}
  Given an undirected, unvaluated graph on $ n $ vertices,
  $ G  = ( N, E ) $, a \textbf{vertex separator set} is a subset
  $ M $ of $ N $, such that removal of the vertices in $ M $
  from $ G $, together with all edges in $ E $ containing them,
  will disconnect the graph.
\end{definition}

For the purposes of this report, consider only graphs that are
connected.  We are interested in vertex separator sets that disconnect
the graph into $2$ subgraphs with approximately equal numbers of
vertices.  We are particularly concerned with \textbf{sparse} graphs and
\textbf{small} separator sets; but we want algorithms that will
automatically perform adequate partitioning of \textbf{any} graph.

Ideally, a partitioning algorithm should find a separator set $ S $
of absolute minimal cardinality, but this appears to be a
combinatorially explosive problem (probably NP-complete), and is
regarded as infeasible.  Instead, we attempt to find an $ S $ that is
``reasonably'' small, at least such that $ | S | \; \ll \; n $,
and partitions $ N $ into two sets of about the same size.

\pagebreak

It is difficult to quantify good separators \cite{Liu:89}.  The ideal
choice is determined by the particular problem, and it is always
possible that a ``small'' separator set does not \textbf{exist} The
extreme example is the completely connected graph, where a separator
set is of order of the same size as the original set.  For planar
graphs, the minimum size of a separator set is always $ \leqslant
\sqrt{8n} $, and the resulting partition has $2$ sides each with no
more than $2/3$ of the total number of vertices \cite{Liu:89}.  This
fact may be of some interest, as many of the large practical problems
are either planar, or very nearly so
(see \cite{DuffGrimesLewis:89,Liu:89,PothenSimonLiou:90}).  Liu,
1989 \cite{Liu:89} describes the notion of ``good'' separator sets in
some detail, and Pothen et al. \cite{PothenSimonLiou:90} provide
several numerical bounds on the minimum possible sizes of separator
sets for a given graph.  Many of the bounds may be grossly
overestimated, and, as they involve calculation of several of the
eigenvalues of the graph Laplacian (an expensive and error-prone
procedure), are not particularly useful.

Methods presented in the literature can be found in the references
to
\cite{%
  DuffGrimesLewis:89,%
  KernighanLin:70,%
  LeisersonLewis:87,%
  Liu:89,%
  PothenSimonLiou:90%
},
and are not reiterated here.  This section deals with a new process,
attributed to \cite{PothenSimonLiou:90}.  The formal construction of
this process depends on theory attributed to \cite{Fiedler:75}, which
begins with a number of relevant definitions.

\begin{definition}
  $
    \mathbf{diameter}
    ( G )
    =
    {\displaystyle \max_{i, j \in N}
    \left\lbrace
      \mathrm{minimum~path~length~} (i, j)
    \right\rbrace}
  $.
\end{definition}

That is, the diameter of a graph $ G = (N, E) $ is the maximum value of
the minimum length path between any two vertices.

\begin{definition}
  A \textbf{decomposition} of a graph $ G = (N, E) $ is a disjoint
  vertex cover $ N = N_1 \cup N_2 $, such that
  $ N_1, N_2 \neq \emptyset $ and $ N_1 \cap N_2 = \emptyset $.
\end{definition}

\begin{definition}
  {$ A $} is called \textbf{irreducible} if, for no decomposition of
  $ N  = N_1 \cup N_2, \; a_{ij} = 0, \; i \in N_1, j \in N_2 $.
  \label{def:Irreducible}
\end{definition}

This corresponds to connectedness of graphs. Using
Definition \ref{def:Irreducible} for symmetric
$ A \in \mathbb{R}^{n \times n} $, define:

\begin{definition}
  $ A $ is of \textbf{degree of reducibility} $ d \in [ 0, n-1 ] $
  if there exists a decomposition of $ N $ into $ d+1 $
  non-empty disjoint subsets
  $
    {\displaystyle
      N = \bigcup_{i = 1}^{d+1} N_i
    }
  $ such that:
  \begin{enumerate}
    \item
      $ A(N_i) $ are irreducible, $ i  = 1, \dots, d+1 $.
    \item
      $ a_{pq} = 0, \; p \in N_i, q \in N_j, i \neq j $.
  \end{enumerate}
\end{definition}

An irreducible symmetric matrix has $ d = 0 $. Thus, a graph has
$ d = 0 $ if it is connected.

\begin{definition}
  The $ n $ eigenvalues of $ A $ (some are possibly multiple) are
  ordered in increasing size:
  $
    \lambda_1
    \leqslant
    \lambda_2
    \leqslant
    \dots
    \leqslant
    \lambda_n
  $.
\end{definition}

Similarly, the corresponding eigenvectors are ordered -- that is, the
$ i $th eigenvector of $ A $ is the eigenvector corresponding to
the $ i $th smallest eigenvalue.  For the purposes of the Graph
Partitioning theorem, it will turn out that all eigenvalues are
non-negative.  The smallest will always be 0, of multiplicity equal to
the number of connected components of the graph (or degree of
reducibility).  For 5-point grid graphs, the largest eigenvalue, by
Theorem \ref{thm:Gershgorin} \cite[p341]{GolubVanLoan:89},
will always be $ \leqslant 8 $.

\begin{theorem}[Gershgorin Circle Theorem]
  Given a matrix $ X^{-1} A X \in \mathbb{C}^{n \times n} $, such that
  $ X^{-1} A X = \mathrm{diag~}(d_1, \dots, d_n ) + F $, where $ F $
  has zero diagonal entries, then:
  \begin{eqnarray*}
    \lambda (A) \subset {\displaystyle \bigcup_{i=1}^n D_i},
    \mathrm{~~where~~}
    D_i =
    \left\lbrace
      z \in \mathbb{C} \mathrm{~such~that~}
      \| z - d_i \| \leqslant
      {\displaystyle \sum_{j=1}^n | f_{ij} |}
    \right\rbrace.
  \end{eqnarray*}
  \label{thm:Gershgorin}
\end{theorem}

That is, the eigenvalues of a complex matrix lie within the union of a
set of $ n $ closed circles in the complex plane.  The centre of each
circle is at the point corresponding to a diagonal element, and its
radius is the sum of the absolute values of the non-diagonal elements
of that row.  Note that it is usual to use $ X = I $.  As applied to
$ A $, the real graph Laplacian of a grid graph, where:

\begin{enumerate}
\item
  The sums of the absolute values of the non-diagonal entries are
  exactly equal to the diagonal elements.
\item
  The diagonal elements of $ A $ are all of value either $2$, $3$ or
  $4$.
\item
  All the eigenvalues of $ A $ are real (as $ A $ is real and
  symmetric).
\end{enumerate}

we obtain $ \lambda (A) \subset [ 0, 8 ] $. More generally, if
$ A $ is a graph Laplacian, $ \lambda(A) \subset [0, 2 \Delta ] $, where
$ \Delta = {\displaystyle \max_i d_i} $.

From Theorems \ref{thm:ZeroEval} and \ref{thm:ZeroEvalplus1}:

\begin{theorem}
  If $ A $ is symmetric, irreducible, has non-negative off-diagonal
  entries, and $ A \mathbf{z} = \mathbf{0} $ for some real
  $n$-vector $\mathbf{z}$, which is neither zero, positive nor
  negative; then $ A $ is not positive semidefinite.
\end{theorem}

Theorem \ref{thm:Corr} is a corollary to a more general result
in \cite{Fiedler:75}.

\begin{theorem}
  If $ A $ is non-negative, irreducible and symmetric, and
  $ \mathbf{y}_d \in \mathbb{R}^n $ is the $ d $th eigenvector of
  $ A $, $ d \geqslant 2 $; then
  $ M = \left\lbrace i \in N | y_i \geqslant 0 \right\rbrace $ is
  non-null and the degree of reducibility of $ A(M) \leqslant d - 2 $.
  \label{thm:Corr}
\end{theorem}

This means that, choosing $ d = 2 $ and $ \mathbf{v} $ as the second
eigenvector of $ A $; the degree of reducibility of
$ A(M) \leqslant 0 $, which means that it \textbf{is} $ 0 $. Using the
second eigenvector of $ A $ as $ \mathbf{v} $ thus yields an irreducible
(connected) component.

Fiedler's paper \cite{Fiedler:75} goes beyond the needs of this report
in defining the graph Laplacian and associated results for the case of
valuated (but not directed) graphs.  The following graph-theoretic
results are a simplification of those presented in Section 3
of \cite{Fiedler:75}.

\begin{definition}
  A \textbf{cut} $ C $ of a graph $ G $ is a set of edges to which a
  decomposition $ N = N_1 \cup N_2 $, where
  $ N_1 \cap N_2 = \emptyset $, exists, such that $ C $ consists
  exactly of all edges in $ G $ with one vertex in each set of the
  decomposition. The subgraphs of $ G $ induced by the subsets
  $ N_1 $ and $ N_2 $ are called \textbf{banks}.
\end{definition}

\begin{theorem}[Unique Decomposition of Connected Banks]
  If there is a decomposition $ N = N_1 \cup N_2 $ of a graph $ G $,
  corresponding to a cut $ C $, such that both corresponding banks are
  connected, then the decomposition of $ N $ corresponding to $ C $
  is unique.
  \label{thm:UniqueConnBanks}
\end{theorem}

\begin{definition}
  The \textbf{algebraic connectivity} \cite{Fiedler:73} of $ L(G) $
  is defined as the smallest non-zero eigenvalue of $ G $.
  Corresponding to this is the \textbf{characteristic valuation},
  which is an assignment of the elements of the eigenvector
  corresponding to the this eigenvalue of $ L(G) $.
\end{definition}

As mentioned in \S\ref{sec:Graph Laplacian}, the smallest
eigenvalue is always $ 0 $, of multiplicity equal to the number of
components (connected units) of $ G $. For connected graphs, the
algebraic connectivity and characteristic valuation are equal to the
second eigenvalue $ \lambda_2 $ and eigenvector $ {\mathbf{v}}_2 $,
respectively.

\begin{theorem}
  For any real $ r $, define
  $ M(r) = \left\lbrace i \in N | y_i \geqslant - r \right\rbrace $.
  The subgraph $ G(r) $ induced by $ G $ on $ M(r) $ is connected.
\end{theorem}

Theorem \ref{thm:MainGraphPartitioning} is the fundamental result
required by the graph partitioning algorithm.

\begin{theorem}[Main Graph-Partitioning Theorem]
  If there exists a real $ c $ such that
  $
    0 \leqslant c < {\displaystyle \max_i y_i}
  $ and
  $
    c \neq y_i \; \forall i
  $, then the set of edges
  $
    (i, j)
  $
  of
  $
    G
  $ for which
  $
    y_i < c < y_j
  $ forms a cut $ C $ of $ G $.
  If
  $
    N_1 = \left\lbrace j \in N | y_j > c \right\rbrace
  $ and
  $
    N_2 = \left\lbrace j \in N | y_j < c \right\rbrace
  $, then
  $N=N_1 \cup N_2$ is a decomposition of
  $N$ corresponding to
  $C$, and the bank
  $G(N_2)$ is connected.
  \label{thm:MainGraphPartitioning}
\end{theorem}

Theorem \ref{thm:UniqueConnBanks} shows that the decomposition and the
banks are \textbf{uniquely} determined. Define
$
  N_1 = \left\lbrace i \in N | y_i > 0 \right\rbrace
$ and
$
  N_2 = \left\lbrace i \in N | y_i < 0 \right\rbrace
$, then $
  N = N_1 \cup N_2
$ is the decomposition corresponding to $ C $.
Theorem \ref{thm:All Cuts Found}, shows that \textbf{all} cuts with
connected banks in a connected graph are able to be obtained (via the
second eigenvector of the graph Laplacian).

\begin{theorem}
  If $ G $ is a connected graph with a cut $ C $ such that both banks
  of $ C $ are connected then there is a positive valuation of the
  edges of $ G $ such that the corresponding characteristic valuation
  $ \mathbf{y} $ is unique (up to a factor), and
  $ y_i \neq 0 \; \forall i $, and $ C $ is formed exactly by
  alternating edges of the valuation $ \mathbf{y} $.
  \label{thm:All Cuts Found}
\end{theorem}

In summary, the essential results in \cite{Fiedler:75} are contained in
Algorithm \ref{alg:EdgeSeparator} (that I have composed), which
is directly referred to in the beginning of
\S\ref{sec:Spectral Partitioning Algorithm}.

\pagebreak

\begin{algorithm}[Edge Separator Algorithm]
  Given an graph $ G = (N, E) $, calculate a edge separator
  set (cut) $ E_1 \subset E $ such that the $ 2 $ resulting
  banks are connected and have approximately equal numbers of vertices.
  \begin{enumerate}
  \item
    Calculate the graph Laplacian $ L(G) $ of $ G $,
    and the smallest non-zero eigenvalue of $ G $ (the
    algebraic connectivity). If $ G $ is connected, the
    eigenvalue $ 0 $ will be of multiplicity $ 1 $, in
    which case, the algebraic connectivity is equal to
    $ \lambda_2 $, the second eigenvalue of $ L(G) $.
  \item
    Calculate the corresponding (second) eigenvector
    (the characteristic valuation).
  \item
    Assign to the vertices of $ G $ the $ n $ elements of the
    characteristic valuation.
  \item
    Find the set of edges $ E_1 \subset E $, whose
    characteristic valuations cross the median value of the
    components of the second eigenvector.
  \item
    $ E_1 $ is the cut required edge separator set.
    Choice of $ E_1 $ from edges whose vertices cross
    some point between two other components of the
    characteristic valuation will also yield $ 2 $
    connected banks, but their sizes will not be nearly equal.
  \end{enumerate}
  \label{alg:EdgeSeparator}
\end{algorithm}

%%%%%%%%%%%%%%%%%%%%%%%%%%%%%%%%%%%%%%%%%%%%%%%%%%%%%%%%%%%%%%%%%%%%%%%%

\subsection{Spectral Partitioning Algorithm}
\label{sec:Spectral Partitioning Algorithm}

The idea of using the results contained in Algorithm
\ref{alg:EdgeSeparator} in an algorithm to partition graphs into vertex
sets appears in Pothen et al., 1990 \cite{PothenSimonLiou:90}. They
describe their algorithm as a ``Spectral Partitioning Algorithm'', and
this convention is followed here. They compare the performance of this
algorithm with several other algorithms:

\begin{enumerate}
\item
  Kernighan-Lin Algorithm -- a modified level structure algorithm
  that is implemented in \textsc{sparspak}%
  \footnote{%
    {\protect\textsc{sparspak}} is a package of sparse matrix routines,
    available through the netlib electronic software library
  }.
\item
  Fiduccia-Mattheyes Algorithm, implemented by Leiserson and
  Lewis \cite{LeisersonLewis:87}.
\item
  Separator Algorithm of Liu \cite{Liu:89}, based on the Multiple
  Minimum Degree Algorithm.
\end{enumerate}

This report is an emulation of their algorithm, formally described in
Algorithm \ref{alg:Spectral Partitioning}, largely taken
from \cite{PothenSimonLiou:90}. This procedure partitions the set of
vertices of a (sparse) graph, in the form specified in
\S\ref{sec:Introduction}. Firstly, it applies the results in
Algorithm \ref{alg:EdgeSeparator} to the graph concerned, to
yield $ E_1 $, an appropriate edge separator set, and
$ A^{\prime} $ and $ B^{\prime} $, sets of vertices on either side of
this set.  Secondly, a combinatorial procedure chooses from the vertices
adjacent to this edge separator set, a vertex-separator set $ S $, and
defines the corresponding vertex sets of the banks $ A $ and $ B $.

Before listing the actual algorithm, the definition of a minimum cover
is required.

\begin{definition}
  Given a graph $ G = (N, E) $, a (vertex) \textbf{cover} is a set
  of vertices $ S $, such that every edge in $ E $ has at
  least one of its endpoints in $ S $.
\end{definition}

\begin{definition}
  A \textbf{minimum cover} is a cover of minimum cardinality.
  \label{def:Minimum Cover}
\end{definition}

Associated with the notion of minimum cover is that of a maximum
matching, which requires another definition.

\begin{definition}
  A \textbf{matching} is a subset of $ E $, such that no two
  endpoints in this subset have the same vertex.
\end{definition}

\begin{definition}
  A \textbf{maximum matching} is a matching of maximal cardinality.
\end{definition}

Maximum matchings and minimum covers are dual concepts, and this is
further discussed in \S\ref{sec:Minimum Cover}.

\begin{figure}
  \begin{algorithm}[Spectral Partitioning Algorithm]
    Given the sparse matrix of a graph $ G = (N, E) $,
    find a partition of the vertex set
    $ N = A \; \cup \; B \; \cup \; S $, such that
    $
      | A | \; \approx \;
      | B | \; \approx \;
      n / 2
    $, and $ | S | $ is ``small'', in a restricted sense.
    \begin{enumerate}
    \item
      Compute the eigenvector $ \mathbf{x}_2 $ and the
      median value $ x_m $ of its components.
    \item
      Partition the vertices of $ G $ into $ 2 $
      sets,
      $
        A^{\prime}
        =
        \left\lbrace v \in N \; | \; x_v \leqslant x_m \right\rbrace
      $ and $ B^{\prime} = N \setminus A^{\prime} $.
      If
      $
        | A^{\prime} | - | B^{\prime} |
        \; > 1
      $, move enough vertices with components equal to
      $ x_m $ from $ A^{\prime} $ to $ B^{\prime} $ to
      make this difference at most one.
    \item
      Define $ A_1 $ as the set of vertices in
      $ A^{\prime} $ adjacent to some vertex in
      $ B^{\prime} $, and $ B_1 $ as the set of
      vertices in $ B^{\prime} $ adjacent to some
      vertex in $ A^{\prime} $. Compute
      $ H = ( A_1, B_1, E_1 ) $, the bipartite subgraph
      induced by the vertex sets $ A_1 $ and $ B_1 $.
    \item
      The required vertex separator set $ S $ is a
      minimum cover of $ H $. It separates $ G $
      into subgraphs with vertex sets $ A =
      A^{\prime} \setminus S $ and $ B = B^{\prime}
      \setminus S $.
    \end{enumerate}
    \label{alg:Spectral Partitioning}
  \end{algorithm}
\end{figure}

Several notes arise in regard to this algorithm:
\begin{enumerate}
\item
  It must be recognised that the ideal aim is to find an $ S $ of
  absolute minimal cardinality, however this algorithm only finds $ S $
  as small as possible in the context of the given edge separator set
  $ E_1 $. In general, the problem of finding the smallest possible
  $ S $ appears to be a combinatorially-explosive one that is not
  achievable by any algorithm efficient enough to be worth considering.
  In consolation, \cite{PothenSimonLiou:90} demonstrate that the
  Spectral Partitioning algorithm generally finds a smaller $ S $ than
  its competitors.
\item
  The problem of finding a minimum cover is non-trivial, and much
  research into it has been performed.
  Pothen et al. \cite{PothenSimonLiou:90} use a ``maximum matching''
  technique (see \S\ref{sec:Minimum Cover}), that is
  \textbf{guaranteed} to give the minimum cover for the given edge
  separator set. The actual implementation in this report uses a
  heuristic procedure to calculate an \textbf{approximate} minimum
  cover, described in Algorithm \ref{alg:AppMinCover} in
  \S\ref{sec:Minimum Cover}. This procedure is \textbf{not}
  guaranteed to give a minimum cover.
\item
  Algorithm \ref{alg:Spectral Partitioning} requires the use of a
  function to determine the median of a list of numbers. An algorithm to
  do this efficiently is a special case of an algorithm to select the
  $ k $th smallest component of a list of numbers. Page 129
  of \cite{Sedgewick:88} describes an algorithm to do this in
  $ O(n \log n ) $ time, by recursively partitioning the list.
\item
  The problem of calculating the second eigenvector of the graph
  Laplacian is also non-trivial, and generally computationally
  expensive. It depends on the prior calculation of the second
  eigenvalue. The most efficient algorithm to find extremal (largest
  and smallest) eigenvalues of sparse matrices is the Lanczos
  algorithm. Unfortunately, the procedure is subject to severe
  numerical problems that make implementation complex, but in practice
  it is the only real choice.  Implementation of the Lanczos algorithm
  is discussed in more detail in \S\ref{sec:LanczosAlgorithm}.
\end{enumerate}

%%%%%%%%%%%%%%%%%%%%%%%%%%%%%%%%%%%%%%%%%%%%%%%%%%%%%%%%%%%%%%%%%%%%%%%%

\subsection{Minimum Cover}
\label{sec:Minimum Cover}

This section deals with the problem of finding a minimum cover (see
Definition \ref{def:Minimum Cover}) of a \textbf{bipartite} graph.
This problem occurs as a necessary component of the Spectral
Partitioning algorithm (Algorithm \ref{alg:Spectral Partitioning}) --
it is desired to find a minimum cover of $ H = (A_1, B_1, E_1 ) $.
One method of doing this is mentioned in \cite{PothenSimonLiou:90}, but
this has not been implemented due to time constraints, and instead a
heuristic procedure is followed. \S\ref{sec:TrueMC} discusses
references in which are found minimum cover algorithms (for both
bipartite and general graphs), and \S\ref{sec:ApproxMC} describes the
heuristic procedure (for bipartite graphs only) actually implemented by
me.

%%%%%%%%%%%%%%%%%%%%%%%%%%%%%%%%%%%%%%%%%%%%%%%%%%%%%%%%%%%%%%%%%%%%%%%%

\subsubsection{True Solution}
\label{sec:TrueMC}

The minimum cover problem has been solved in a number of ways. The
earliest algorithms for the minimum cover of a bipartite graph are
based on the ``Dulmage-Mendelsohn
decomposition'' \cite{DulmageMendelsohn:63,JohnsonDulmageMendelsohn:62},
but these references are not particularly readable. A more general
result, for the minimum cover of any graph is found
in \cite{NormanRabin:59}, but for the purposes of this partitioning, it
is better to only consider algorithms for bipartite graphs, to maximise
efficiency. For the rest of this section, consider the term graph to
mean the bipartite graph $ H = (A_1, B_1, E_1) $. Also define
$ a = \; | A_1 | $, $ b = \; | B_1 | $,
$ e = \; | E_1 | $ and $ n = a + b $.

As mentioned in \S\ref{sec:Spectral Partitioning Algorithm},
the minimum cover of a bipartite graph is the dual of the maximum
matching, although the details of this relationship are not explicitly
supplied. Pothen et al. \cite{PothenSimonLiou:90} cite a further
paper \cite{PothenFan:90} that details an algorithm for a maximum
matching. A simplified description of this algorithm is found in
\cite[pp221--227]{PapadimitriouSteiglitz:82}, and is reported to solve
the matching problem in $ O(\min(a, b).e) $ time.

An alternative approach is presented in pages 495--499
of \cite{Sedgewick:88}, and deals with the matching problem in terms of
a problem in ``flow maximisation''. The algorithm described is quoted
as requiring $ O( n (e + n) \log n) $ time (or $ O(n^3) $ time, for
dense graphs). It is expected that $ H $ is in general sparse, but
this algorithm does not appear to be as efficient as that used
in \cite{PothenSimonLiou:90}.

Implementation of a true minimum cover algorithm is left as a future
exercise, see \S\ref{sec:Further Work}.

\pagebreak

%%%%%%%%%%%%%%%%%%%%%%%%%%%%%%%%%%%%%%%%%%%%%%%%%%%%%%%%%%%%%%%%%%%%%%%%

\subsubsection{Heuristic Solution}
\label{sec:ApproxMC}

The heuristic procedure which I have used is formally described in
Algorithm \ref{alg:AppMinCover}.

Input data is an edge separator set $ E_1 $, generated by the
vertices whose eigenvaluation crosses the median eigenvalue. Construct
$ d $, a listing of the degrees of the vertices in
$ H = A_1 \cup B_1 $, and then, whilst any element in $ d $ is
positive, perform the following procedure: Find the current vertex of
maximum degree in $ d $, add it to $ S $, and remove it and the edges
adjacent to it (in $ E_1$), from $ d $. Repeat until no edges of degree
$ > 0 $ remain in $ d $. Setting $ d_j $ to $ 0 $ ensures that the
vertex will not be considered again. Implementation is improved by the
very cheap process of comparing the size of the resultant $ S $ with
the sizes of $ A_1 $ and $ B_1, $ and if $ S $ is larger than the
smaller of these, $ S $ is replaced by the smaller one. Since no
elements are removed from $ S $, it can be implemented as a simple list
or array.

\begin{algorithm}[Approximate Minimum Cover Algorithm]
  Given $ E_1 $, the\\
  edge separator set of edges of a bipartite
  graph on vertex sets $ A_1 $ and $ B_1 $, find a vertex
  (separator) set $ S $ such that each edge in $ E_1 $ is
  incident to a vertex in $ S $, with $ | S | $ ``small''.
  \begin{tabbing}
    ~~~~~~~~~ \=~~~~~~~~~ \= \kill
    $ S \leftarrow \emptyset; $ \\
    $ \textbf{for~} i \in A_1 \cup B_1 $ \\
      \> $ d_i \leftarrow \mathrm{degree~of~} i
        \mathrm{~in~} E_1; $ \\
    $ \textbf{while~} \mathrm{some~} d_i > 0 $ \\
      \> $ \mathrm{choose~} j \mathrm{such~that~}
        d_j = {\displaystyle \max_i d_i}; $ \\
      \> $ d_j \leftarrow 0; $ \\
      \> $ S \leftarrow S \cup \left\lbrace j \right\rbrace; $ \\
      \> $ \textbf{for~} \mathrm{each~} i \mathrm{adjacent~to~} j
        \mathrm{~in~} E_1 $ \\
      \>  \> $ d_i \leftarrow d_i - 1; $ \\
    $ \mathrm{if~} | S | \; > \; | A_1 | $ \\
      \> $ S \leftarrow A_1; $ \\
    $ \mathrm{if~} | S | \; > \; | B_1 | $ \\
      \> $ S \leftarrow B_1; $
  \end{tabbing}
  \label{alg:AppMinCover}
\end{algorithm}

This strategy will work to varying degrees with different examples. In
practice, it would be possible to create a (possibly pathological)
example for which this algorithm would create an unreasonably large
$ S $. For sparse graphs of large diameter, it is expected that
$ | A_1 | \; \ll \; | A^{\prime} | $ and
$ | B_1 | \; \ll \; | B^{\prime} | $, so the process for
finding $ S $ should yield the desired set $ | S | \; \ll \; n $.
It is \textbf{guaranteed} only that
$
    | S | \; \leqslant \;
    \min ( | A_1 |, | B_1 | )
$.

\pagebreak

%%%%%%%%%%%%%%%%%%%%%%%%%%%%%%%%%%%%%%%%%%%%%%%%%%%%%%%%%%%%%%%%%%%%%%%%

\subsection{Lanczos Algorithm}
\label{sec:LanczosAlgorithm}

The Lanczos algorithm (originally attributable to \cite{Lanczos:50}) is
an efficient method of finding the extremal eigenvalues of a sparse
matrix. As the Spectral Partitioning algorithm is to be applied to
sparse matrices, the Lanczos algorithm is required in its
implementation. It is listed formally in Algorithm \ref{alg:Lanczos},
copied almost verbatim from Chapter 9 of \cite{GolubVanLoan:89}.

\begin{figure}
  \begin{algorithm}[The Lanczos Algorithm]
    Given a symmetric $ A \in \mathbb{R}^{n \times n} $ and
    $ \mathbf{w} \in \mathbb{R}^n $ having unit 2-norm,
    compute a symmetric, tridiagonal matrix
    $ T_j \in \mathbb{R}^{j \times j} $ with the property that
    $ \lambda ( T_j ) \subset \lambda(A) $. The diagonal and
    subdiagonal elements of $ T_j $ are stored in
    $ \alpha ( 1:j ) $ and $ \beta (1:j-1) $ respectively.

    \begin{tabbing}
    ~~~~~~~~~ \=~~~~~~~~~ \= \kill
    $
      \mathbf{v} = \mathbf{0}  \; \; \;
      {\beta}_0 = 1; \; \; \; j = 0;
    $
    \\
    $
      \textbf{while~} {\beta}_j \neq 0
    $
    \\
      \> $ \textbf{if~} j \neq 0 $ \\
      \>  \>
        $
          \mathbf{t} = \mathbf{w}; \; \; \;
          \mathbf{w} = {\displaystyle \frac{1}{{\beta}_j}} \mathbf{v};
          \; \; \;
          \mathbf{v} = - {\beta}_j \mathbf{t};
        $
        \\
      \> $ \textbf{end} $ \\
      \>
      $
        \mathbf{v} = \mathbf{v} + A \mathbf{w}; \; \; \;
        j = j + 1;
      $
      \\
      \>
      $
        {\alpha}_j = {\mathbf{w}}^{\top} \mathbf{v}; \; \; \;
        \mathbf{v} = \mathbf{v} - {\alpha}_j \mathbf{w}; \; \; \;
        {\beta}_j = {\| \mathbf{v} \|}_2;
      $
      \\
      $ \textbf{end} $
    \end{tabbing}
    \label{alg:Lanczos}
  \end{algorithm}
\end{figure}

Practical implementation of the Lanczos algorithm almost always
requires some reorthogonalisation. For the purposes of the Spectral
Partitioning algorithm, we need only reorthogonalise against
$\mathbf{e}$ (the $ n $-vector of ones), as the subspace required must
be orthogonal to this. Rounding error can be avoided by the following
strategy: When computing $ \mathbf{y} = A \mathbf{x} $, instead of
simply returning $\mathbf{y}$, we will return $\mathbf{y}$
orthogonalised against $\mathbf{e}$:

\begin{eqnarray*}
  \mathbf{y} \leftarrow
  \mathbf{y} -
  {\displaystyle \frac{
    \mathbf{e}^{\top}\mathbf{y}
  }{
    n
  }}
  \mathbf{e}.
\end{eqnarray*}

That is, return
$
  ( I - \mathbf{e}^{\top}\mathbf{e}/n )\/ A \mathbf{x}
$ instead of $ A \mathbf{x} $, as well as starting with
$ \mathbf{x}_0 \perp \mathbf{e} $. The latter requirement is satisfied
by the choice of
$ {(\mathbf{x}_0)}_i = i - (n + 1)/2, \; \; i = 1, \dots, n $, a choice
recommended in \cite{PothenSimonLiou:90}.

\pagebreak

%%%%%%%%%%%%%%%%%%%%%%%%%%%%%%%%%%%%%%%%%%%%%%%%%%%%%%%%%%%%%%%%%%%%%%%%

\subsection{Recursive Decomposition}

It is possible, and may be necessary for many applications, to be able
to repeat the decomposition process on the subcomponents $ A $ and
$ B $ of $ N $. Once the Spectral Partitioning algorithm is
implemented, it can be recursively called, until the subgraphs $ A $
and $ B $ are of a certain manageable ``atomic'' size
(about $3$ -- $10$ vertices), yielding a complete permutation of the
vertices in $ N $, via the procedure described in Algorithm
\ref{alg:Recursive}.

\begin{algorithm}[Recursive Decomposition Algorithm]
  Given a graph $ G = (N, E) $, on the set of vertices
  $N=\left\lbrace 1, \dots, n \right\rbrace $, decompose $ N $ into a
  permutation of itself, such that the smallest size units in the
  permutation are no larger than $ p $, the ``atomic'' size.
  \begin{tabbing}
  ~~~~~~~~~ \=~~~~~~~~~ \= \kill
  $ \textbf{function~} partition(N) $ \\
    \> Apply the Spectral Partitioning Algorithm
      to partition $ N $ into $ A \cup B \cup S; $ \\
    \> $ \textbf{if~} | A | \; > \; p $ \\
    \>   \> $ A \leftarrow partition(A); $ \\
    \> $ \textbf{if~} | B | \; > \; p $ \\
    \>   \> $ B \leftarrow partition(B); $ \\
    \> $ \textbf{return~} ( A \; | \; B \; | \; S ) $
  \end{tabbing}
  \label{alg:Recursive}
\end{algorithm}

\pagebreak

%%%%%%%%%%%%%%%%%%%%%%%%%%%%%%%%%%%%%%%%%%%%%%%%%%%%%%%%%%%%%%%%%%%%%%%%

\section{Implementation}
\label{sec:Implementation}

The algorithms presented in \S\ref{sec:Theory} are implemented as a small
collection of functions and programs in ANSI-standard \textsc{C} that are
interfaced with a \textsc{C} software library (\textsc{meschach} \cite{Stewart:91}).
They were written within a \textsc{unix} (BSD4.3) environment.
Appendix \ref{app:Code} provides a listing of the source code, although the
functions and data-structures in \textsc{meschach} that are referenced are
not explicitly listed, and the reader is referred
to \cite{Stewart:91}.

Sections \ref{sec:decomp} to \ref{sec:testdc} describe the
operation and application of each of the programs/functions listed in
Appendix \ref{app:Code}. Currently, there is no documentation on
these codes apart from this section and comments in the relevant source
code. \S\ref{sec:Results and Problems} discusses the results
of applying these functions to the sample graphs described in
Appendix \ref{app:Graphs}.

%%%%%%%%%%%%%%%%%%%%%%%%%%%%%%%%%%%%%%%%%%%%%%%%%%%%%%%%%%%%%%%%%%%%%%%%

\subsection{The Function \texttt{decomp}}
\label{sec:decomp}

The function \texttt{decomp} is the primary partitioning (and recursive
decomposition) routine written. Input/output parameters are:

\begin{enumerate}
\item
  \texttt{sp\_mat *L}: Pointer to the sparse matrix of the graph
  Laplacian of the graph $ G $ that we wish to partition, of
  size $ n $. Unchanged on exit.
\item
  \texttt{PERM *P}: On entry, a pointer to a permutation (list) of
  length $ n $ of the actual numbers of the vertices in the set
  being partitioned. $ P $ is returned as the permutation of
  itself corresponding to the partition.
\item
  \texttt{PERM *A, PERM *B, PERM *S}: The components of the
  partition. Currently not actually relevant, but will be used in
  future developments, where the structure of a recursive
  decomposition will also be returned. These are pointers to
  permutations of size $ n $ on entry, containing dummy data.
  Returned as correctly-sized and filled permutations.
\item
  \texttt{int rec\_lvl}: \texttt{decomp} is designed perform one of two
  types of tasks:
  \begin{enumerate}
  \item
    Partition a sparse graph once.
  \item
    Recursively repeat this until the units involved are of
    size less than a nominated (currently hard-wired into
    \texttt{decomp} as 3) ``atomic'' size. $ P $ is returned
    as the permutation which will be most useful for
    subsequent factorisation. This is not directly useful,
    at present, as there is no record of the \textbf{structure}
    of $ P $ -- again, this is left as future
    work, and is discussed in \S\ref{sec:Further Work}.
  \end{enumerate}
  \texttt{rec\_lvl} is used to tell \texttt{decomp} which of these 2
  tasks to perform. If \texttt{rec\_lvl} is set to -1 by a driver
  program, \texttt{decomp} will only partition the graph once, any
  other choice allowing it to recursively partition the set until
  satisfied. The parameter is used as a record of the depth of
  recursion by \texttt{decomp}, so a driver program will typically
  only ever set \texttt{rec\_lvl} to $0$ or $-1$.
\end{enumerate}

Currently, \texttt{decomp} is called by the driver program
\texttt{testdc} (see \S\ref{sec:testdc}), and the parameter
\texttt{rec\_lvl} is set by \texttt{testdc} according to one of its
input (\texttt{argv}) arguments.

%%%%%%%%%%%%%%%%%%%%%%%%%%%%%%%%%%%%%%%%%%%%%%%%%%%%%%%%%%%%%%%%%%%%%%%%

\subsection{The Function \texttt{mk\_sp\_graph}}
\label{sec:mk_sp_graph}

A function to generate random sparse graphs. It is not usually used for
testing \texttt{decomp}, as the graphs generated are possibly not
connected, and are not of small diameter, which means that we cannot
expect small separator sets from them. Future versions of this function
should generate essentially planar graphs of large diameter, which we
could expect to partition into sets with small $ S $. Input parameters
are:

\begin{enumerate}
\item
  \texttt{unsigned int n}: The order of the graph desired.
\item
  \texttt{unsigned int p}: The average degree of the vertices.
\end{enumerate}

The function returns a pointer to a \texttt{sp\_mat}, and is called by
\texttt{testdc}, a driver program for \texttt{decomp}.

%%%%%%%%%%%%%%%%%%%%%%%%%%%%%%%%%%%%%%%%%%%%%%%%%%%%%%%%%%%%%%%%%%%%%%%%

\subsection{The Function \texttt{select}}
\label{sec:select}

This function selects the \texttt{k}th smallest element of a series of
(real) numbers, stored in a \texttt{VEC *}, (a pointer to a
\texttt{VEC}). It is an implementation of an efficient algorithm
in \cite[page 128]{Sedgewick:88}. In particular, it can be used to
determine the median element of a series of numbers.

Input parameters are:
\begin{enumerate}
\item
  \texttt{VEC *y}: The vector involved.
\item
  \texttt{int k}: \texttt{select} returns the \textbf{value} (not the
  position) of the \texttt{k}th smallest element of \texttt{y}.
\end{enumerate}

\texttt{select} returns a \texttt{double}. \texttt{decomp} calls it to
find the second smallest element in the vector of ``eigenvalues''
returned by \texttt{trieig}. The function is also called by
\texttt{decomp} to determine the second smallest eigenvalue in a list
(see Algorithm \ref{alg:Spectral Partitioning}).

%%%%%%%%%%%%%%%%%%%%%%%%%%%%%%%%%%%%%%%%%%%%%%%%%%%%%%%%%%%%%%%%%%%%%%%%

\subsection{The Program \texttt{gr\_lap}}
\label{sec:gr_lap}

This program generates sparse matrices of graphs associated with
solving Laplace's Equation on a rectangular 5-point grid. The graphs,
not the graph Laplacians are generated. \texttt{gr\_lap} automatically
creates a file of the appropriate name, which can be changed for later
use. For example (in a \textsc{unix} environment), typing the command:

\texttt{~~~~~~~~gr\_lap 13 23}\\
creates a file called ``\texttt{Lap.13.23}'' in the current directory,
which contains (in standard format for \textsc{meschach} to read), the
sparse graph related to solving Laplace's Equation on a
$ 13 \times 23 $ grid.

\pagebreak

%%%%%%%%%%%%%%%%%%%%%%%%%%%%%%%%%%%%%%%%%%%%%%%%%%%%%%%%%%%%%%%%%%%%%%%%

\subsection{The Program \texttt{testdc}}
\label{sec:testdc}

This program is the main driver written for \texttt{decomp}, and is
called with a series of (\texttt{argv}) input parameters. For example
(in a \textsc{unix} environment), typing the command:

\texttt{~~~~~~~~testdc 157 7 13 -1}\\
calls \texttt{testdc}, and tells it that we wish to use a graph on $157$
vertices, with an average vertex degree of $7$, and to only print out
intermediate results of data-structures of size $13$ or less.
\texttt{testdc} will then prompt the user to select either one of a
range of graphs stored in the current directory, or generate a random
sparse graph, for partitioning. If the user chooses to input a graph
that does not have the correct dimension, \texttt{testdc} will abort.
This stupid-seeming system is also intended to allow the user to
instruct \texttt{testdc} to generate random sparse graphs (using
\texttt{mk\_sp\_graph}), with parameters \texttt{n} $= 157$, and
\texttt{p} $= 7$. Thus, the user is expected to have some knowledge of
the database of sparse graphs before using \texttt{testdc}. The system
will be improved in subsequent versions of \texttt{testdc} and
\texttt{decomp}.

The final parameter is the value of \texttt{rec\_lvl} (see
\S\ref{sec:decomp}) that we wish to initially pass to \texttt{decomp}.
If set to -1, \texttt{decomp} is instructed to only partition the graph
once. Any other value, or omission of it instructs \texttt{decomp} to
recursively decompose the graph into atomic subunits. The atom-size is
currently hard-wired into \texttt{decomp} (as 3), but could become an
input parameter in future versions. \texttt{testdc} currently calls
\texttt{select} and \texttt{decomp}, as well as numerous subroutines
from \textsc{meschach}.

\pagebreak

%%%%%%%%%%%%%%%%%%%%%%%%%%%%%%%%%%%%%%%%%%%%%%%%%%%%%%%%%%%%%%%%%%%%%%%%

\section{Results and Problems}
\label{sec:Results and Problems}

This section details some results obtained by

\begin{enumerate}
\item
  Partitioning of all the graphs listed in Table \ref{tab:GraphListing}
  (page \pageref{tab:GraphListing} in Appendix \ref{app:Graphs}).
\item
  Full recursive decomposition of the $5$ smallest ones.
\end{enumerate}

It refers to the example graphs listed in Appendix \ref{app:Graphs}.
These graphs were used in the development of the programs, and have been
used to illustrate of the progress of the algorithms. There is only one
picture of a grid graph -- for Ishmail, on $55$ vertices. The other grid
graphs are too large to depict. The successful partitioning of these
graphs is mentioned in \S\ref{sec:Partitioning Results},
recursive decomposition in \S\ref{sec:Recursion Results}, and
problems encountered are discussed in \S\ref{sec:Problems}.

Direct comparison with other published results (in
particular \cite{PothenSimonLiou:90}) is not possible, as the computers
involved are of very different speeds. The test problems that have been
used are largely from the Boeing-Harwell sparse matrix test problem
library \cite{DuffGrimesLewis:89}. Access paths to this database were
discovered too late for it to be used. Other implementational problems,
such as the Boeing-Harwell database being stored as a column-oriented
data-structure, made the application of the partitioning algorithm to
this suite currently infeasible. (\textsc{meschach} deals primarily in
terms of row-oriented data-structures.)

%%%%%%%%%%%%%%%%%%%%%%%%%%%%%%%%%%%%%%%%%%%%%%%%%%%%%%%%%%%%%%%%%%%%%%%%

\subsection{Partitioning}
\label{sec:Partitioning Results}

\texttt{decomp}, driven by \texttt{testdc}, has been used to correctly
partition all of the graphs listed in Table \ref{tab:GraphListing}.
The times (in CPU seconds) taken to do this on the University of
Queensland's Mathematics Department Pyramid 9810 computer (operating
under \textsc{unix} BSD4.3) are listed in Table \ref{tab:CPUtimes}. The
Pyramid 9810 is benchmarked by a \textsc{linpack}%
\footnote{%
  {\protect\textsc{linpack}} is available through the netlib electronic
  software library
}
routine at approximately $0.5$ Mflop s$^{-1}$. The computer used by
Pothen et al. \cite{PothenSimonLiou:90} is a
{\protect\textsc{cray y-mp}},%
\footnote{%
  Trademark of {\protect\textsc{cray}} Research
}
and is expected to be benchmarked about $2$ orders of magnitude faster
than the Pyramid 9810. Results are roughly comparable with those
in \cite{PothenSimonLiou:90}, by incorporating a scaling factor of
about $100$.

All but one of the separator sets listed in Table \ref{tab:CPUtimes}
are of the minimum possible size. The separator set of size $101$ for
Schlomo could have been improved (to size 61), by changing tolerances
used in \texttt{decomp}, but this would have required prohibitive
amounts of memory.

\begin{table}
  \centering
  \begin{tabular}{||l|r|r|l|r||} \hline\hline
    Name & \protect$ n \protect$ &
    \protect$ | S | \protect$ &
    \protect$ S \protect$ & CPU time (s) \\ \hline\hline
    Idit & 7 & 1 & 3 & 0.0 \\
    Moshe & 20 & 2 & 9, 10 & 0.1 \\
    Itzchak & 31 & 1 & 30 & 0.2 \\
    Shimuel & 32 & 2 & 7, 11 & 0.1 \\
    Ishmail & 55 & 5 & 5, 16, 27, 38, 49 & 0.1 \\
    Yacov & 105 & 5 & 9, 31, 52, 73, 94 & 0.3 \\
    Yair & 121 & 11 & 59--69 & 0.2 \\
    Arieh & 505 & 5 & 50, 151, 252, 353, 454 & 8.1 \\
    Aaron & 2121 & 21 & 50, 151, 252, \dots, 2070 & 22.0 \\
    Schlomo & 6161 & 101 & 3079--3179 & 44.4 \\
    Shimshon & 6400 & 80 & 3120--3199 & 57.8 \\ \hline\hline
  \end{tabular}
  \caption{
    Times (in CPU seconds) and \protect$ | S | \protect$
    for {\protect\tt decomp} to partition graphs listed in
    Table \protect\ref{tab:GraphListing}.
  }
  \label{tab:CPUtimes}
\end{table}

%%%%%%%%%%%%%%%%%%%%%%%%%%%%%%%%%%%%%%%%%%%%%%%%%%%%%%%%%%%%%%%%%%%%%%%%

\subsection{Full Recursive Decomposition}
\label{sec:Recursion Results}

The full recursive decomposition was successfully completed for the
graphs Idit, Moshe, Itzchak, Shimuel, Ishmail, and Yacov, finding
permutations of the vertices that were traced to be exactly what would
be expected for the first four cases. \texttt{decomp} was allowed to
run on Ishmail and Yacov until it had completed a full recursive
decomposition -- this ran to $6$ levels of recursion for Yacov. Results
appeared correct, although were not manually checked explicitly! The
full recursion on the other (very large) problems was not attempted.

Caveat: the above comments refer to a \textbf{previous} version of {\tt
decomp}, not the one supplied in Appendix \ref{app:Code}, which
has some minor implementational problems. Currently, \texttt{decomp}
crashes at some point during recursive calls.

%%%%%%%%%%%%%%%%%%%%%%%%%%%%%%%%%%%%%%%%%%%%%%%%%%%%%%%%%%%%%%%%%%%%%%%%

\subsection{Problems}
\label{sec:Problems}

This section describes the most significant problem encountered in
implementing the Spectral Partitioning algorithm
(Algorithm \ref{alg:Spectral Partitioning} on
page \pageref{alg:Spectral Partitioning}). This problem is in the
choice of $ m $, the order of the tridiagonal matrix $ T_m $
returned by the Lanczos algorithm.

The first impression is that any choice of $ m $ will do, the larger
the better, in estimating $ \lambda_2 $ accurately. The constraint is
that the computational expense is dependent on $ m^2 $ (amongst other
things), and $ m $ cannot be increased without limit. At the very
least, $ m \leqslant n $ would appear to be an obvious, if
exaggerated bound. There is then, an optimal choice of $ m $, such that:

\begin{enumerate}
\item
  $ \lambda_2 $ is found accurately.
\item
  Minimal computational expense is involved.
\end{enumerate}

In fact, empirical observation demonstrates that there is typically a
range of suitable $ m $, this range varying in position and bandwidth
with the problem encountered. It is \textbf{not} possible to accurately
predict this range in advance, although typically it is located in the
vicinity of $ \sqrt{n} $. ($ n $ is the number of vertices in the
graph.)

To make matters more difficult, if $ m $ is chosen from outside this
range, the computed value of $ \lambda_2 $ will be \textbf{wrong}. It
appears that if $ m $ is too small, $ \lambda_2 $ will be too
large, sometimes not even corresponding to larger eigenvalues of the
graph Laplacian. If $ m $ is too large, the calculated $ \lambda_2 $
cannot possibly be an eigenvalue (as it is supposed to be the second
smallest one!), and with increasing $ m $ the returned values of
$ \lambda_2 $ typically reduce, suddenly hitting 0, and remaining there.
This phenomenon is called the ``Ghost Eigenvalue''
problem \cite{GolubVanLoan:89}, and has been studied but apparently not
yet conquered.

Thus, a strategy is needed to choose an optimal $ m $. Empirical
observations suggest that if $ m $ is increased from below the
optimal position, when the resulting computed values of $ \lambda_2 $
converge, they converge to the \textbf{correct} value. This is not
necessarily true when decreasing $ m $ from above the optimum, but it
may be so.

Several approaches are possible in developing a strategy for choosing
$ m$:

\begin{enumerate}
\item
  Begin with a value of $ m $ \textbf{much} larger than is
  expected to be adequate for the given problem. For example, if
  $ m = 20 $ is guessed to be optimum, try $ m = 40 $. Using
  this value, run the Lanczos algorithm and investigate the
  resulting $ \lambda_2 $. As the Lanczos algorithm has
  returned us $ T_m $, by decrementing $ m $ and considering
  $ T(1:m-1,1:m-1) $, it is possible to calculate
  $ \lambda_2 $ using $ m-1 $. This procedure is repeated
  until the computed $ \lambda_2 $ converge. It is not expensive,
  as no new calculations of $ T_m $ have to be made. The problems
  with this method are that:

  \begin{enumerate}
  \item
    It is not certain that convergence from above generates
    the correct $ E_1 $.
  \item
    Calculation of $ T_m $ for an overly large $ m $ is
    expensive.
  \item
    The choice of the initial $ m $ is problem-dependent.
    If this choice cannot be automated, there is no hope of
    an automatic sparse matrix factorisation technique
    being derived. The only obvious way to automate this
    choice would be to use $ m = n $, and this would be
    so expensive as to make the whole technique infeasible.
  \end{enumerate}

  This strategy has been considered, but was not implemented due
  to these concerns. Its prime advantage is in its ease of
  programming.
\item
  The second method is to begin with a small $ m $, say
  $ m = \min ( n, 3 ) $, and always maintain $ m $ at least this
  size. Run the Lanczos algorithm using this value of $ m $,
  and return. First calculate the resulting $ \lambda_2 $, and
  then the value that $ \lambda_2 $ would have been if $ m $
  were one less. If these two values are the same, accept $ m $
  and $ \lambda_2 $, else increment $ m $ by $2$ and re-run the
  Lanczos algorithm. This technique is suited to calling the
  Lanczos algorithm as an external function, but suffers from
  being wasteful in its computational requirements, as the $ T_m $
  calculations have to be repeated each time $ m $ is
  incremented.  In practice, this technique has been experimented
  with, and works, but has since been superseded by the next technique,
  which is much cheaper to implement, although more difficult to
  program.
\item
  Follow a similar strategy to the above, but maintain all the
  information on the Lanczos algorithm as it progresses. In order
  to do this, a modified Lanczos algorithm is incorporated into
  the part of the code that examines the convergence of the
  computed second eigenvalue, so that the need for external function
  calls is removed. This technique is incorporated into the
  current version of \texttt{decomp}, that appears in
  Appendix \ref{app:Code}.
\end{enumerate}

If the problem being considered is on $ n $ vertices, with average
degree $ k $, and $ m $ is the actual size of the computed $ T $
matrix, then the computational cost of the second and third methods can
be compared by the following analysis. \cite{GolubVanLoan:89} describes
the cost of one Lanczos iteration (without reorthogonalisation) as
$(2k+8)n$ flops.

For the second method, returning $ T_j $ requires $ ( 2 k + 8 ) n j $
flops. If this cost is repeated for $ j = 3, 5, \dots, m $, then
the total cost is:

\begin{eqnarray*}
  {\displaystyle
    \sum_{j = 3}^m
    ( 2 k + 8 ) n j
  }
  =
  {\displaystyle
    \sum_{i = 1}^{(m-1)/2}
    ( 2 k + 8 ) n ( 2 i + 1)
  }
  \approx
  ( k + 4 ) n m^2 / 2 \mathrm{~flops}.
\end{eqnarray*}

For example, for the $ 21 \times 101 $ $5$-point grid graph Aaron,
where $ k \approx 4, n = 2121 $ and observation suggests that
$ m = 35 $ is a good choice, this represents $ \approx 10^7 $ flops.

Compare this with the cost for the third method, which is simply the
cost of the last Lanczos algorithm call of the second (plus some small
overhead). This is approximately
$ ( 2 k + 8 ) m n \approx 2 \times 10^6 $ flops, an improvement of
order $ m / 4 $. This means, for example, using the Pyramid 9810
computer, running at $ \approx 0.5 \times 10^6 $ flops s$ ^{-1} $,
calculating the second eigenvector takes $2$ rather than $20$ CPU
seconds.

In summary, none of the methods described are \textbf{guaranteed} to
work, as they are all based on the observation that when the computed
second eigenvalues converge, they converge to the correct value of
$\lambda_2 $. Thus their philosophy is heuristic, and needs refinement,
but it \textbf{appears} to work for the examples tested. This is an
outstanding problem for further work.

\pagebreak

%%%%%%%%%%%%%%%%%%%%%%%%%%%%%%%%%%%%%%%%%%%%%%%%%%%%%%%%%%%%%%%%%%%%%%%%

\section{Further Work}
\label{sec:Further Work}

This section is a listing of directions for further work, and includes
some commentary of the work of \cite{PothenSimonLiou:90}:

The incorporation of the Spectral Partitioning algorithm into a
complete function to factorise and solve large sparse linear systems
using nested dissection is a primary research target, and would depend
on a number of subsidiary goals, which are described below.

\begin{enumerate}
\item
  Problems in dealing with the choice of $ m $ in the Lanczos
  algorithm have been discussed in
  \S\ref{sec:Problems}. Research is required to
  establish a more rigorous solution to this problem.
\item
  \texttt{decomp} should be augmented to become a function that also
  returns the \textbf{structure} of the permutation $ P $. This
  might be made possible by using a ternary code associated with
  each element of $ P $, to describe the level of recursion
  required to obtain it. This would facilitate a factorisation
  routine that recursively decomposed the large sparse matrix
  into manageable atomic-sized units. Implementation would be
  coupled with an analysis of the performance of this technique
  relative to other methods, such as Gauss-Seidel and Successive
  Over Relaxation.
\item
  Currently, there is currently very little error-checking in any
  of the functions and programs that I have written, but this
  facility is readily implemented in the context of the interface
  with \textsc{meschach}. Input parameters are not carefully checked
  for existence and correctness of dimension, and there is not
  much checking through the codes for other problems, especially
  numerical problems. In part, the latter is due to time
  constraints, but it is also due to the total complexity of the
  algorithm. Further work is needed to make the routines more
  robust.
\item
  More efficient data-structures would improve programming and
  execution speed. In particular, the edge listing $ H $, is
  currently represented as an $ | H | \times \; 2 $ array
  of \texttt{double}s (\texttt{MAT *H}). (The same is true for the
  degree listing $ D $.) As they are intended to be integers,
  this means that input, output and comparison operations on
  their elements require repeated casts. This is poor programming, and
  consumes more CPU time than required. A replacement system could
  involve simple integer arrays, or more powerfully, it could be a
  \textsc{C} structure containing an array, each element of which
  points to $2$ integers representing adjacent edges, modelled after
  other \textsc{meschach} structures.  A possible data-structure and
  referencing description of this system is provided in Figure
  \ref{fig:EdgeList}.

  In conjunction with this we could develop a graph-input
  routine, together with a suite of error-checking routines.
  Currently, the graph is input as an adjacency matrix. This
  could be exchanged for an incidence matrix input. \texttt{decomp}
  could internally convert this to the appropriate adjacency
  matrix. This would improve the ability of the user to input a
  graph accurately. The driver program could be extended to
  handle a series of different input formats.

\pagebreak

  Other possibilities include storing the adjacency matrix of the
  graph as a dense \textbf{bit} array, rather than a sparse integer
  \textsc{C} structure. This would however, make implementation far more
  difficult, and is not considered to be practical for the
  internal workings of \texttt{decomp}. The programs \texttt{gr\_lap}
  and \texttt{testdc}, which involve input and output of graphs (as
  their adjacency matrices), could be made considerably more
  efficient by this practice.

\begin{figure}
  \centering
  {\tt
  \begin{verbatim}
        typedef struct edge_elt
        {
                unsigned int    first, second;
        } edge_elt

        typedef struct edge_list
        {
                unsigned int    m, max_m;
                ee              *edge_elt;
        } edge_list;

       PERM *decomp( ... )
        {
                edge_list       *D, *H;
        /* "(H->ee[i]).first" now replaces "H->me[i][0]" */
  \end{verbatim}
  }
  \caption{
     {\textsc{C}}
     Data structure for edge-listing and degree-recording.
  }
  \label{fig:EdgeList}
\end{figure}

\item
  The nested dissection algorithm lends itself to implementation
  on parallel processing machines, and a long-term goal could be
  to implement it in such an environment.
\item
  Pothen et al. \cite{PothenSimonLiou:90} discuss various
  estimates for the size of an adequate separator set, calculable
  in terms of various eigenvalues of the graph Laplacian and
  other parameters of the graph. It might be valuable to
  investigate the use of these estimates for some internal
  consistency checks within the decomposition algorithm.
\item
  The implementation of a true minimum cover algorithm, as
  described in \S\ref{sec:TrueMC}, would ensure that
  $ | S | $ is as small as possible. The algorithm presented
  in \cite[pp221--227]{PapadimitriouSteiglitz:82} appears
  to be the most efficient and accessible choice. It could be
  expected that the increase in computational expense in using
  this minimum cover algorithm would be compensated for by the
  reduction in computational expense (due to the smaller
  $ | S | $) in a full factorisation routine.
\end{enumerate}

\pagebreak

%%%%%%%%%%%%%%%%%%%%%%%%%%%%%%%%%%%%%%%%%%%%%%%%%%%%%%%%%%%%%%%%%%%%%%%%

\section{Listing of Graphic Examples}
\label{app:Graphs}

This appendix is a listing and discussion of some of the graphic
examples used in developing the code and algorithms. Names and
descriptions of the experimental graphs are presented in
Table \ref{tab:GraphListing}. Should these names appear unfamiliar,
the reader is informed that they are loosely transliterated Hebrew
names of various people in the Hebrew Bible. There is no particular
significance in this choice of names.

\begin{table}[ht]
  \centering
  \begin{tabular}{||l|r|l||} \hline\hline
    Name & \protect$ n \protect$ & Description \\ \hline\hline
    Idit & 7 & Simplest case, almost a tree \\
    Moshe & 20 & Hand-drawn example \\
    Itzchak & 31 & Almost a tree \\
    Shimuel & 32 & Simple example \\
    Ishmail & 55 & \protect$ 5 \times 11 \protect$ 5-point grid \\
    Yacov & 105 & \protect$ 5 \times 21 \protect$ 5-point grid \\
    Yair & 121 & \protect$ 11 \times 11 \protect$ 5-point grid \\
    Arieh & 505 & \protect$ 5 \times 101 \protect$ 5-point grid \\
    Aaron & 2121 &
    \protect$ 21 \times 101 \protect$ 5-point grid \\
    Schlomo & 6161 &
    \protect$ 61 \times 101 \protect$ 5-point grid \\
    Shimshon & 6400 &
    \protect$ 80 \times 80 \protect$ 5-point grid \\ \hline\hline
  \end{tabular}
  \caption{Listing of the experimental graphs.}
  \label{tab:GraphListing}
\end{table}

The large grid graphs Schlomo ($61 \times 101 $) and Shimshon
($80 \times 80 $) are used for comparison with the results
of \cite{PothenSimonLiou:90}.

%%%%%%%%%%%%%%%%%%%%%%%%%%%%%%%%%%%%%%%%%%%%%%%%%%%%%%%%%%%%%%%%%%%%%%%%

\subsection{Idit}

\begin{figure}[htbp]
  \begin{center}
    \epsfig{file=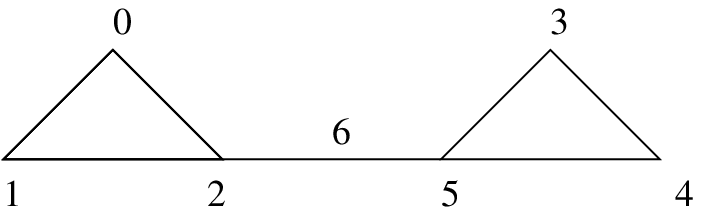}
    \caption{Idit.}
    \label{fig:Idit}
  \end{center}
\end{figure}

\begin{figure}[htbp]
  \begin{center}
    \epsfig{file=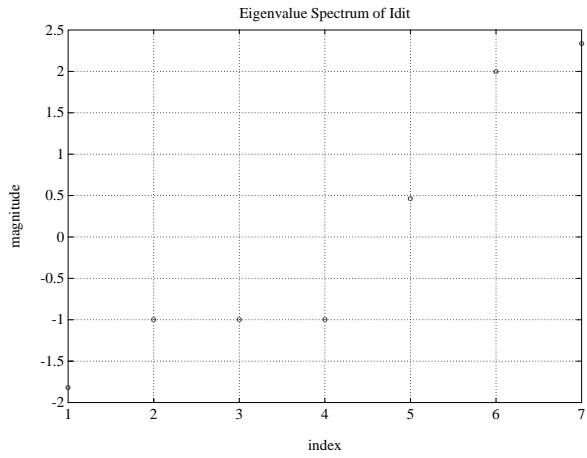,height=175pt}
    \caption{The eigenvalue spectrum of Idit.}
    \label{fig:Idit_eval}
  \end{center}
\end{figure}

\begin{figure}[htbp]
  \begin{center}
    \epsfig{file=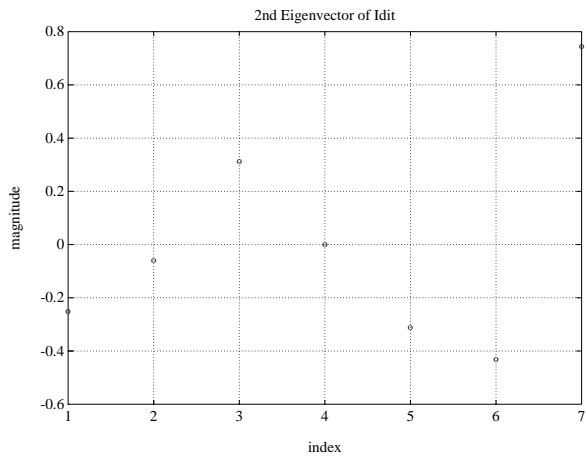,height=175pt}
    \caption{The second eigenvector of Idit.}
    \label{fig:Idit_evec}
  \end{center}
\end{figure}

Idit is the smallest graph examined ($7$ vertices, $8$ edges), and is
depicted in Figure \ref{fig:Idit}. Its graph Laplacian is presented
in Figure \ref{fig:IditGrLap}. Figures \ref{fig:Idit_eval}
and \ref{fig:Idit_evec} are plots of the eigenvalue spectrum and the
second eigenvector (respectively) of the graph Laplacian of Idit, as
generated by \textsc{matlab}.%
\footnote{%
  {\protect\textsc{matlab}} is an (interpreted) matrix computation
  package, and is a trademark of The Mathworks, Inc.
}

\begin{figure}
  \centering
  $
    L =
    \left[
    \begin{array}{*{7}{r}}
       2 & -1 & -1 &  . &  . &  . &  . \\
      -1 &  2 & -1 &  . &  . &  . &  . \\
      -1 & -1 &  3 & -1 &  . &  . &  . \\
       . &  . & -1 &  2 & -1 &  . &  . \\
       . &  . &  . & -1 &  3 & -1 & -1 \\
       . &  . &  . &  . & -1 &  2 & -1 \\
       . &  . &  . &  . & -1 & -1 &  2
      \end{array}
  \right]
  $
  \caption{The sparse graph Laplacian of Idit.}
  \label{fig:IditGrLap}
\end{figure}

\clearpage

%%%%%%%%%%%%%%%%%%%%%%%%%%%%%%%%%%%%%%%%%%%%%%%%%%%%%%%%%%%%%%%%%%%%%%%%

\subsection{Moshe}
\label{sec:Moshe}

\begin{figure}[htbp]
  \begin{center}
    \epsfig{file=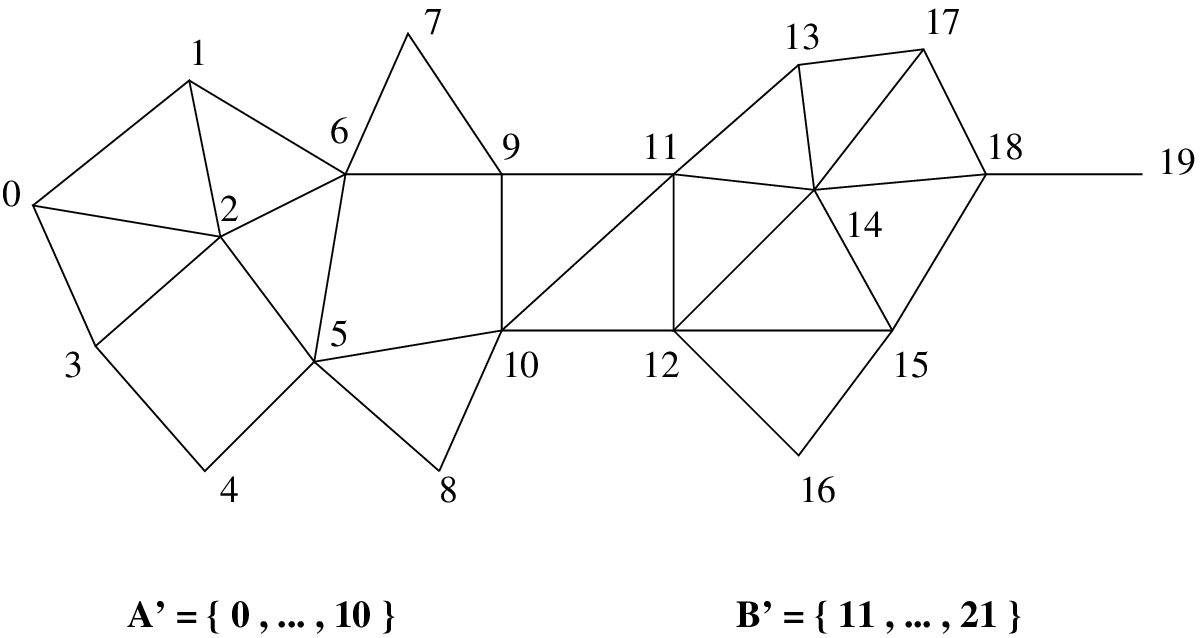,width=180pt}
    \caption{Moshe.}
    \label{fig:Moshe}
  \end{center}
\end{figure}

Moshe is a hand-drawn graph made for illustration of the progress of
the Spectral Partitioning algorithm. It is a planar graph ($20$
vertices, $32$ edges, diameter $7$), and is depicted in Figure
\ref{fig:Moshe}.  The operation of the Spectral Partitioning algorithm
(Algorithm \ref{alg:Spectral Partitioning}) is described for Moshe in
the following paragraphs.

The initial graph partitioning process partitions
$N=\left\lbrace 0,\dots, 19 \right\rbrace $, into $2$ almost equal
halves,  $ A^{\prime} = \left\lbrace 0, \dots, 9 \right\rbrace $ and
$B^{\prime} = \left\lbrace 10, \dots, 19 \right\rbrace $, with an
\textbf{edge} separator set of size 3:
$ E_1 = \left\lbrace (8, 10), (9, 10), (9, 11) \right\rbrace $.  $A_1$
is the set of vertices in $ A^{\prime} $ that are adjacent to vertices
in $ B^{\prime} $, and vice-versa. Inspection shows that the relevant
bipartite graph on $ E_1 $ consists of
$ A_1 = \left\lbrace 8, 9 \right\rbrace $ and
$ B_1 = \left\lbrace 10, 11 \right\rbrace $.

Using only $ A_1 $, $ B_1 $ and $ E_1 $, we seek to calculate a
\textbf{vertex} separator set $ S $, as a subset of the vertices in
$ A_1 $ and $ B_1 $, such that all edges in $ A_1 $ and $ B_1 $
are incident upon at least one vertex in $ S $. We would like $ S $
to be of minimum cardinality. Whilst this may in general be a
non-trivial problem, in this case inspection shows that
$
  \left\lbrace 9, 10 \right\rbrace, \; \left\lbrace 8,
  9 \right\rbrace \mathrm{~and~} \left\lbrace 10, 11 \right\rbrace
$
are minimum covers, the first case corresponding to
Algorithm \ref{alg:AppMinCover} accepting the initial choice of
$ S $, the others corresponding to the use of $ A_1 $ or $ B_1 $
as the cover.

\begin{figure}[htbp]
  \begin{center}
    \epsfig{file=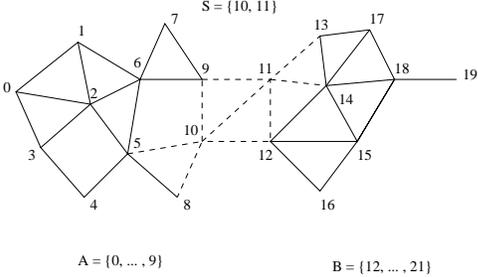,width=180pt}
    \caption{Moshe, after partitioning.}
    \label{fig:Moshe2}
  \end{center}
\end{figure}

Lastly, find $ A = A^{\prime} \setminus S $ and
$ B = B^{\prime} \setminus S $. Use of
$ S = \left\lbrace 9, 10 \right\rbrace $ gives
$ A = \left\lbrace 0, \dots, 8 \right\rbrace $ and
$ B = \left\lbrace 11, \dots, 19 \right\rbrace $. This choice is
depicted in Figure \ref{fig:Moshe2}.  The resultant set of cut edges
are drawn with dashed lines.

\begin{figure}[htbp]
  \begin{center}
    \epsfig{file=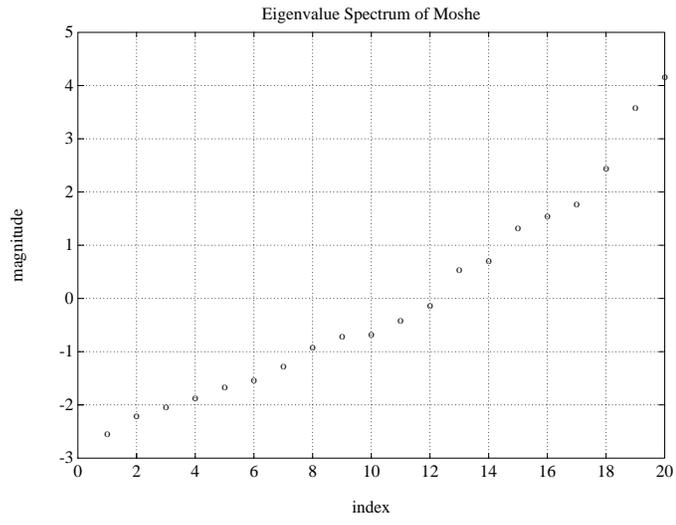,height=200pt}
    \caption{The eigenvalue spectrum of Moshe.}
    \label{fig:Moshe_eval}
  \end{center}
\end{figure}

\begin{figure}[htbp]
  \begin{center}
    \epsfig{file=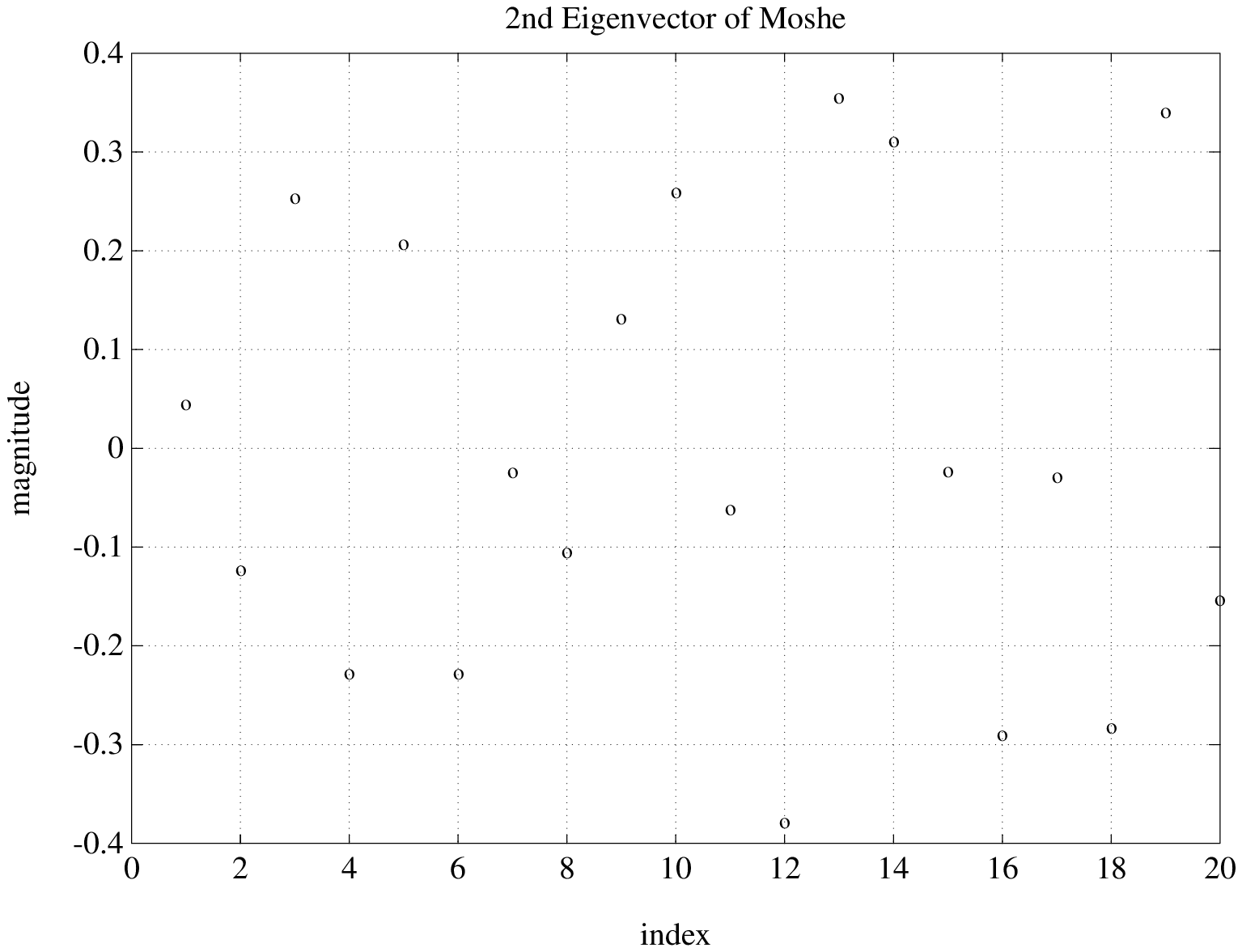,height=200pt}
    \caption{The second eigenvector of Moshe.}
    \label{fig:Moshe_evec}
  \end{center}
\end{figure}

Figures \ref{fig:Moshe_eval} and \ref{fig:Moshe_evec} are plots
of the eigenvalue spectrum and the second eigenvector (respectively) of
the graph Laplacian of Moshe, as generated by \textsc{matlab}.

\clearpage

%%%%%%%%%%%%%%%%%%%%%%%%%%%%%%%%%%%%%%%%%%%%%%%%%%%%%%%%%%%%%%%%%%%%%%%%

\subsection{Itzchak}

\begin{figure}[htbp]
  \begin{center}
    \epsfig{file=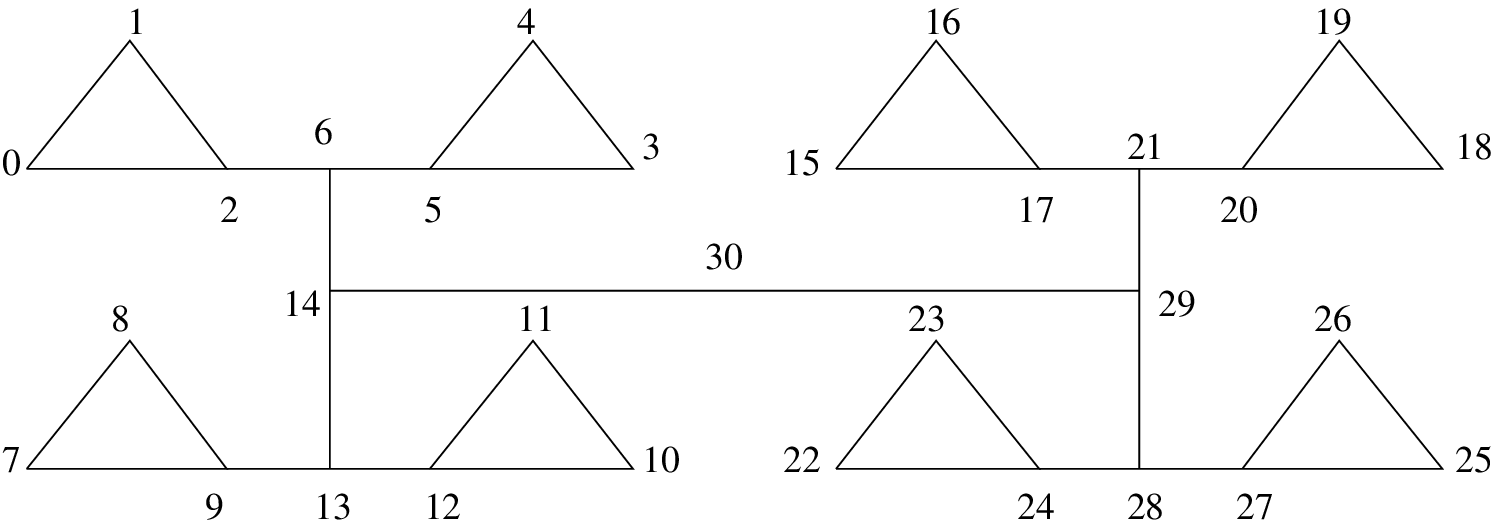,width=360pt}
    \caption{Itzchak.}
    \label{fig:Itzchak}
  \end{center}
\end{figure}

\begin{figure}[htbp]
  \begin{center}
    \epsfig{file=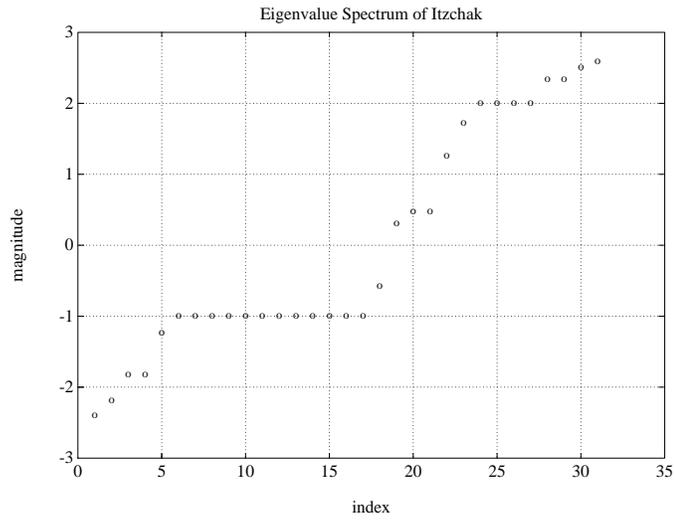,height=200pt}
    \caption{The eigenvalue spectrum of Itzchak.}
    \label{fig:Itzchak_eval}
  \end{center}
\end{figure}

\begin{figure}[htbp]
  \begin{center}
    \epsfig{file=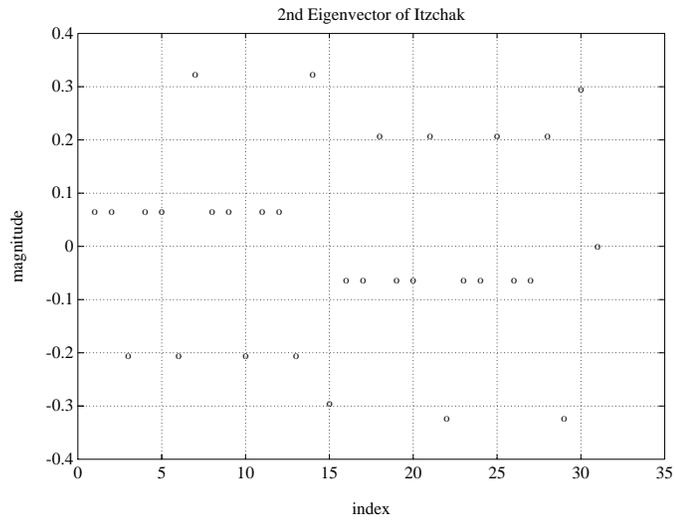,height=200pt}
    \caption{The second eigenvector of Itzchak.}
    \label{fig:Itzchak_evec}
  \end{center}
\end{figure}

The graph Itzchak is presented in Figure \ref{fig:Itzchak}. Its
adjacency matrix is quite illustrative, and is presented in
Figure \ref{fig:ItzchakMatrix}. Figures \ref{fig:Itzchak_eval}
and \ref{fig:Itzchak_evec} are plots of the eigenvalue spectrum and
the second eigenvector (respectively) of the graph Laplacian of Itzchak,
as generated by \textsc{matlab}.

\def\sss{\scriptscriptstyle}

\begin{figure}[ht]
\begin{eqnarray*}
  \hspace{-60pt}
  \left[
    \begin{array}{*{2}{*{2}{*{2}{*{3}{r}|}r|}r|}r}
      &{\sss 1}&{\sss 1}& & & & & & & & & & & & & & & & & & & & & & & & & & & & \\
      {\sss 1}& &{\sss 1}& & & & & & & & & & & & & & & & & & & & & & & & & & & & \\
      {\sss 1}&{\sss 1}& & & & & {\sss 1} & & & & & & & & & & & & & & & & & & & & & & & & \\ \hline
      & & & &{\sss 1}&{\sss 1}& & & & & & & & & & & & & & & & & & & & & & & & & \\
      & & & {\sss 1}& &{\sss 1}& & & & & & & & & & & & & & & & & & & & & & & & & \\
      & & & {\sss 1}&{\sss 1}& & {\sss 1} & & & & & & & & & & & & & & & & & & & & & & & & \\ \hline
      & &{\sss 1}& & &{\sss 1}& & & & & & & & &{\sss 1}& & & & & & & & & & & & & & & & \\ \hline
      & & & & & & & &{\sss 1}&{\sss 1}& & & & & & & & & & & & & & & & & & & & & \\
      & & & & & & &{\sss 1}& &{\sss 1}& & & & & & & & & & & & & & & & & & & & & \\
      & & & & & & &{\sss 1}&{\sss 1}& & & & &{\sss 1}& & & & & & & & & & & & & & & & & \\ \hline
      & & & & & & & & & & &{\sss 1}&{\sss 1}& & & & & & & & & & & & & & & & & & \\
      & & & & & & & & & & {\sss 1}& &{\sss 1}& & & & & & & & & & & & & & & & & & \\
      & & & & & & & & & & {\sss 1}&{\sss 1}& &{\sss 1}& & & & & & & & & & & & & & & & & \\ \hline
      & & & & & & & & &{\sss 1}& & &{\sss 1}& &{\sss 1}& & & & & & & & & & & & & & & & \\ \hline
      & & & & & &{\sss 1}& & & & & & &{\sss 1}& & & & & & & & & & & & & & & & &{\sss 1}\\ \hline
      & & & & & & & & & & & & & & & &{\sss 1}&{\sss 1}& & & & & & & & & & & & & \\
      & & & & & & & & & & & & & & &{\sss 1}& &{\sss 1}& & & & & & & & & & & & & \\
      & & & & & & & & & & & & & & &{\sss 1}&{\sss 1}& & & & &{\sss 1}& & & & & & & & & \\ \hline
      & & & & & & & & & & & & & & & & & & &{\sss 1}&{\sss 1}& & & & & & & & & & \\
      & & & & & & & & & & & & & & & & & &{\sss 1} & &{\sss 1}& & & & & & & & & & \\
      & & & & & & & & & & & & & & & & & &{\sss 1} &{\sss 1}& &{\sss 1}& & & & & & & & & \\ \hline
      & & & & & & & & & & & & & & & & &{\sss 1}& & &{\sss 1}& & & & & & & & &{\sss 1}& \\ \hline
      & & & & & & & & & & & & & & & & & & & & & & &{\sss 1}&{\sss 1}& & & & & & \\
      & & & & & & & & & & & & & & & & & & & & & &{\sss 1}& &{\sss 1}& & & & & & \\
      & & & & & & & & & & & & & & & & & & & & & &{\sss 1}&{\sss 1}& & & & &{\sss 1}& & \\ \hline
      & & & & & & & & & & & & & & & & & & & & & & & & & &{\sss 1}&{\sss 1}& & & \\
      & & & & & & & & & & & & & & & & & & & & & & & & &{\sss 1} & &{\sss 1}& & & \\
      & & & & & & & & & & & & & & & & & & & & & & & & &{\sss 1} &{\sss 1}& &{\sss 1}& & \\ \hline
      & & & & & & & & & & & & & & & & & & & & & & & &{\sss 1}& & &{\sss 1}& &{\sss 1}& \\ \hline
      & & & & & & & & & & & & & & & & & & & & &{\sss 1}&& & & & & &{\sss 1}& &{\sss 1}\\ \hline
      & & & & & & & & & & & & & &{\sss 1}& & & & & & & && & && & & &{\sss 1}&
    \end{array}
  \right]
  \end{eqnarray*}
  \caption{The adjacency matrix $L$ of Itzchak.}
  \label{fig:ItzchakMatrix}
\end{figure}

\pagebreak

%%%%%%%%%%%%%%%%%%%%%%%%%%%%%%%%%%%%%%%%%%%%%%%%%%%%%%%%%%%%%%%%%%%%%%%%

\subsection{Shimuel}

\begin{figure}[htbp]
  \begin{center}
    \epsfig{file=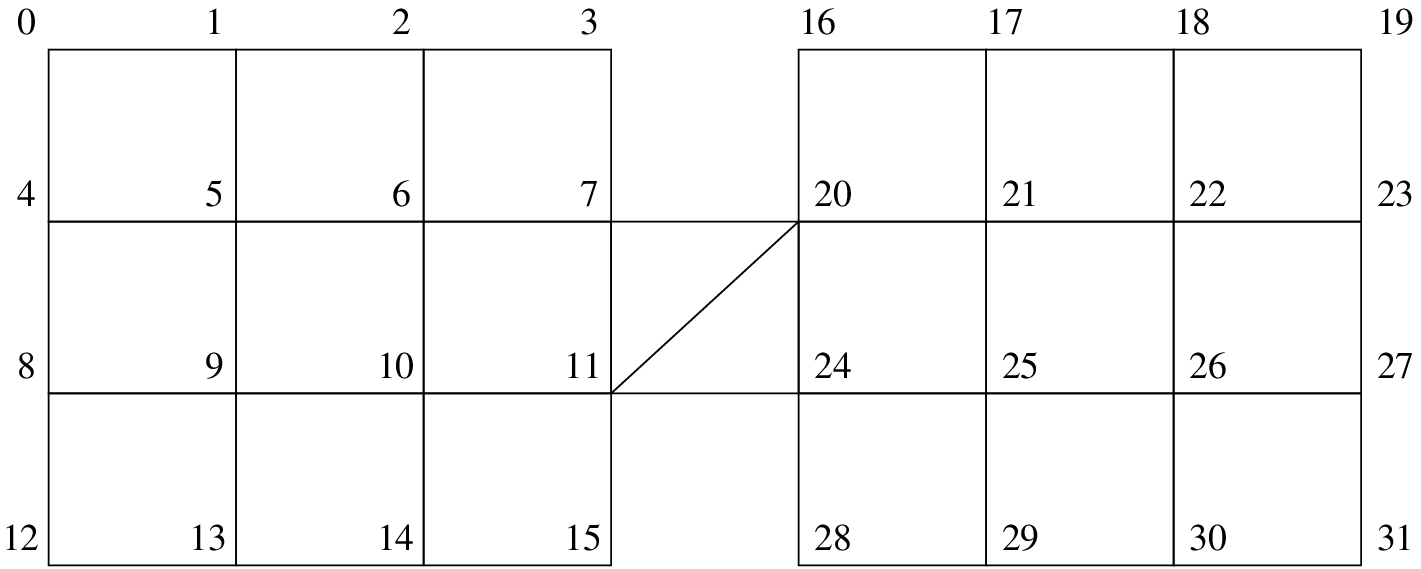,width=360pt}
    \caption{Shimuel.}
    \label{fig:Shimuel}
  \end{center}
\end{figure}

\begin{figure}[htbp]
  \begin{center}
    \epsfig{file=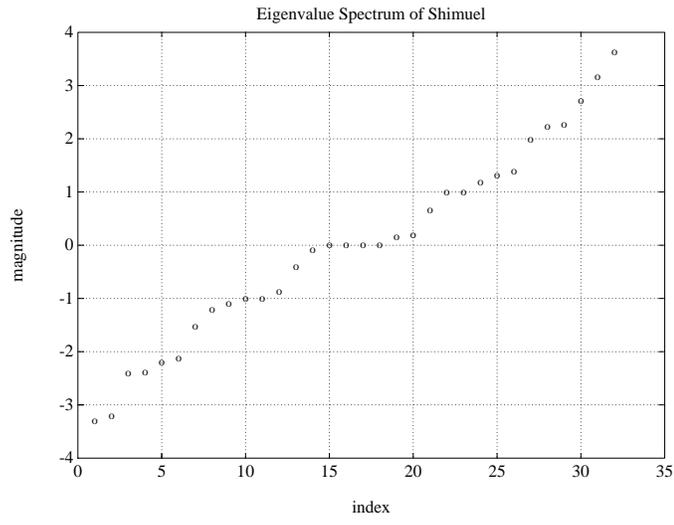,height=200pt}
    \caption{The eigenvalue spectrum of Shimuel.}
    \label{fig:Shimuel_eval}
  \end{center}
\end{figure}

\begin{figure}[htbp]
  \begin{center}
    \epsfig{file=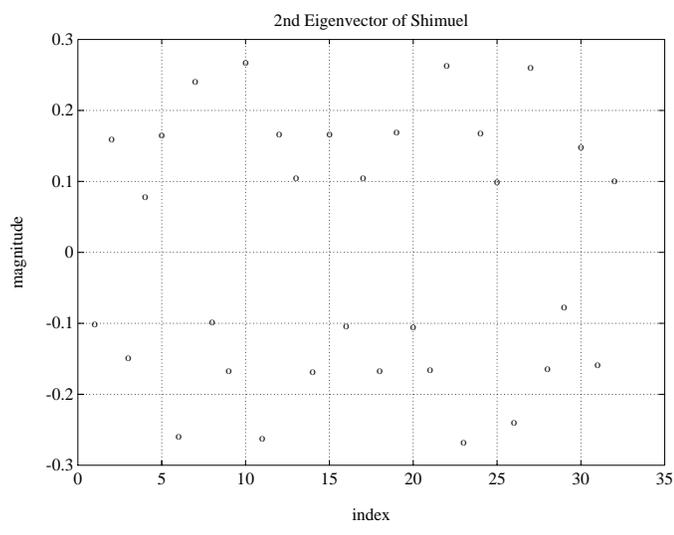,height=200pt}
    \caption{The second eigenvector of Shimuel.}
    \label{fig:Shimuel_evec}
  \end{center}
\end{figure}

Shimuel (Figure \ref{fig:Shimuel}) is another example on a small number
of vertices, also used for development purposes. It is set up to find
the separator set $ S = \left\lbrace 11, 20 \right\rbrace $.
Figures \ref{fig:Shimuel_eval} and \ref{fig:Shimuel_evec} are
plots of the eigenvalue spectrum and the second eigenvector
(respectively) of the graph Laplacian of Shimuel, as generated by
\textsc{matlab}.

\clearpage

%%%%%%%%%%%%%%%%%%%%%%%%%%%%%%%%%%%%%%%%%%%%%%%%%%%%%%%%%%%%%%%%%%%%%%%%

\subsection{Ishmail}
\label{sec:Ishmail}

\begin{figure}[htbp]
  \begin{center}
    \epsfig{file=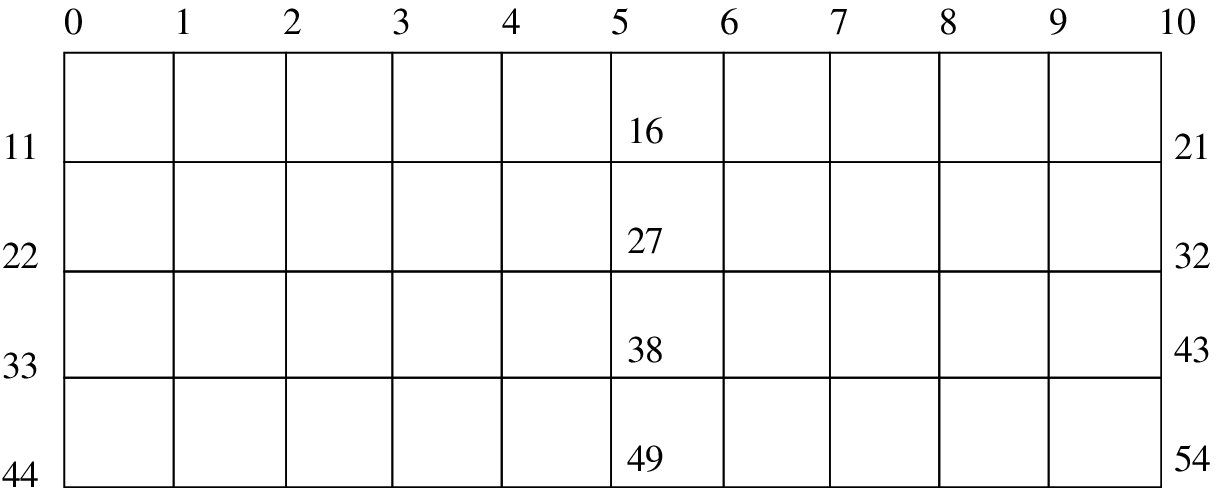,width=270pt}
    \caption{Ishmail.}
    \label{fig:Ishmail}
  \end{center}
\end{figure}

Ishmail is the only (5-point) grid graph represented in this appendix,
and is depicted in Figure \ref{fig:Ishmail}. It is a $ 5 \times 11 $
grid ($ n = 55 $ vertices). The adjacency matrix of Ishmail
illustrates that referred to in \S\ref{sec:Graph Laplacian},
and its form is:

\begin{eqnarray*}
  G
  =
  \left[
    \begin{array}{*{5}{c}}
    P   & I_5 & .   & .   & .   \\
    I_5 & P   & I_5 & .   & .   \\
    .   & I_5 & P   & I_5 & .   \\
    .   & .   & I_5 & P   & I_5 \\
    .   & .   & .   & I_5 & P
    \end{array}
  \right],
\end{eqnarray*}
where
\begin{eqnarray*}
  P
  =
  \left[
    \begin{array}{*{11}{c}}
    . & 1 & . & . & . & . & . & . & . & . & . \\
    1 & . & 1 & . & . & . & . & . & . & . & . \\
    . & 1 & . & 1 & . & . & . & . & . & . & . \\
    . & . & 1 & . & 1 & . & . & . & . & . & . \\
    . & . & . & 1 & . & 1 & . & . & . & . & . \\
    . & . & . & . & 1 & . & 1 & . & . & . & . \\
    . & . & . & . & . & 1 & . & 1 & . & . & . \\
    . & . & . & . & . & . & 1 & . & 1 & . & . \\
    . & . & . & . & . & . & . & 1 & . & 1 & . \\
    . & . & . & . & . & . & . & . & 1 & . & 1 \\
    . & . & . & . & . & . & . & . & . & 1 & .
    \end{array}
  \right].
\end{eqnarray*}

\begin{figure}[htbp]
  \begin{center}
    \epsfig{file=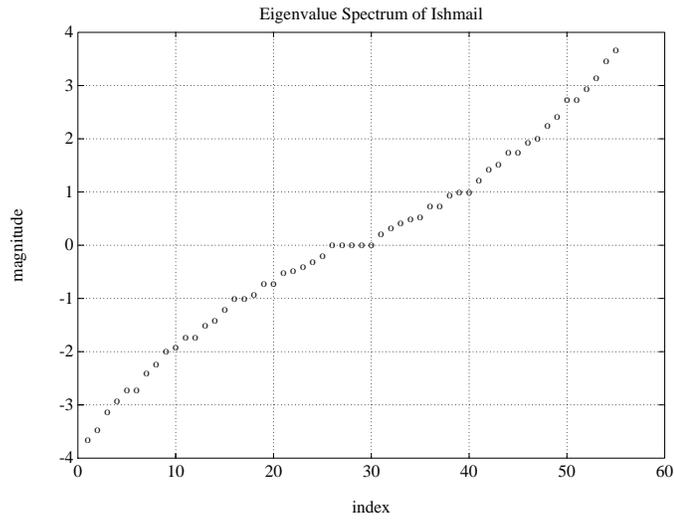,height=200pt}
    \caption{The eigenvalue spectrum of Ishmail.}
    \label{fig:Ishmail_eval}
  \end{center}
\end{figure}

\begin{figure}[htbp]
  \begin{center}
    \epsfig{file=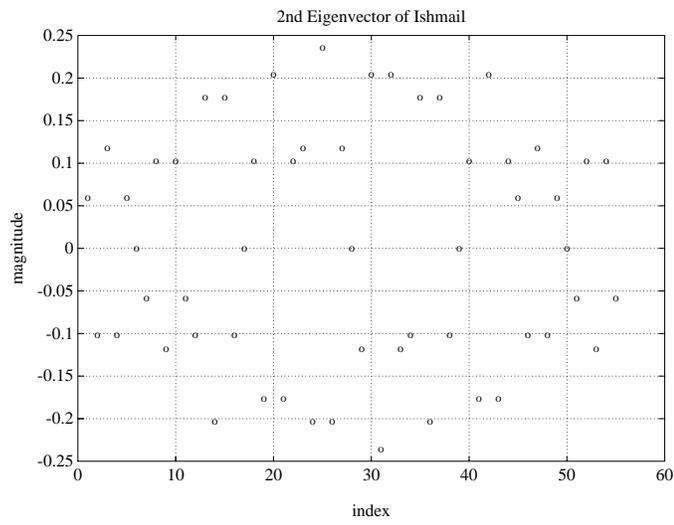,height=200pt}
    \caption{The second eigenvector of Ishmail.}
    \label{fig:Ishmail_evec}
  \end{center}
\end{figure}

Figures \ref{fig:Ishmail_eval} and \ref{fig:Ishmail_evec} are
plots of the eigenvalue spectrum and the second eigenvector
(respectively) of the graph Laplacian of Ishmail, as generated by
\textsc{matlab}.

Section 4 (titled ``Graph Products'') of \cite{PothenSimonLiou:90},
discusses the expected repetition of the second eigenvector of an
$ m \times n $ 5-point grid graph, showing that it can be found as the
Kronecker (tensor, outer) product of the length-$ m $ path graph and
an $ \mathbf{e} $-vector of size $ n $. Examination of
Figure \ref{fig:Ishmail_evec} (and Figure \ref{fig:Yacov_evec})
clearly demonstrates this result. Clever exploitation of this property
may reduce the overall computational expense required for the Spectral
Partitioning algorithm, but not significantly, as the expense is
dominated by the computations that yield $ \lambda_2 $ (see
\S\ref{sec:Problems}). As this property is only true for
\textbf{grid} graphs, implementation would require \texttt{decomp} to
have an extra input flag.

\clearpage

%%%%%%%%%%%%%%%%%%%%%%%%%%%%%%%%%%%%%%%%%%%%%%%%%%%%%%%%%%%%%%%%%%%%%%%%

\subsection{Yacov}

\begin{figure}[htbp]
  \begin{center}
    \epsfig{file=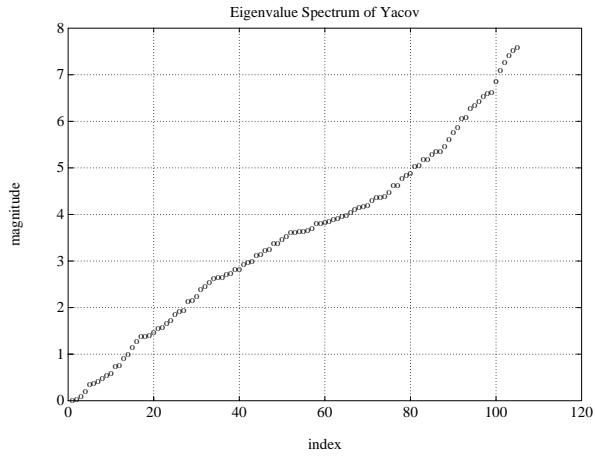,height=175pt}
    \caption{The eigenvalue spectrum of Yacov.}
    \label{fig:Yacov_eval}
  \end{center}
\end{figure}

\begin{figure}[htbp]
  \begin{center}
    \epsfig{file=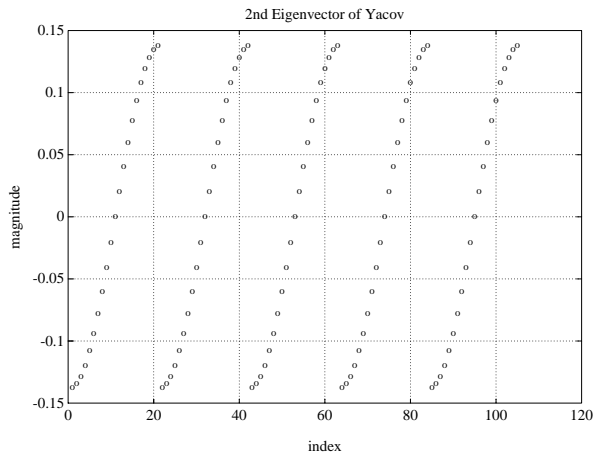,height=175pt}
    \caption{The second eigenvector of Yacov.}
    \label{fig:Yacov_evec}
  \end{center}
\end{figure}

Yacov is a $ 5 \times 21 $ 5-point grid graph on $ n = 105 $ vertices,
and is not depicted here. Figures \ref{fig:Yacov_eval}
and \ref{fig:Yacov_evec} are plots of the eigenvalue spectrum and
the second eigenvector (respectively) of the graph Laplacian of Yacov,
as generated by \textsc{matlab}.

\pagebreak

%%%%%%%%%%%%%%%%%%%%%%%%%%%%%%%%%%%%%%%%%%%%%%%%%%%%%%%%%%%%%%%%%%%%%%%%

\section{Code Listing}
\label{app:Code}

As mentioned in \S\ref{sec:Implementation}, the algorithms presented in
this report are implemented in ANSI-standard \textsc{C}, and interface
with \textsc{meschach}, the \textsc{C} software library for numerical
analysis written by David Stewart \cite{Stewart:91}. The contents of
this appendix are listings of the source code that I have written for
this work. Without access to \cite{Stewart:91}, references to data
structures and functions called from \textsc{meschach} will be
meaningless. Listings of the data-structures and functions are not
provided here as they are quite long. \S\ref{sec:Problems} mentions
that the code written is not necessarily bug-free, and currently has
very little error-handling capacity.

\subsection{Listing of \texttt{decomp}}
\scriptsize\tt 
\begin{verbatim}
/*                            decomp                                  */


/* *********************************************************************


Contents  : This file is called "decomp.c", and contains the function
		"decomp".

     Aim  : (Recursively) decompose the vertex set of a (sparse) graph,
	    represented by a graph Laplacian L, into separator sets.
	    This function returns a pointer to a permutation, which,
	    when applied to L, will decompose it for efficient
	    recursive LDL^T factorisation.

 Language : ANSI Standard C

  Author  : David  De Wit (and David Stewart)  March 5  -  June 19  1991


********************************************************************* */


#include <math.h>
#include <stdio.h>
#include "matrix.h"
#include "matrix2.h"
#include "sparse.h"
#include "sparse2.h"


extern double   select(VEC *, u_int);
extern int      icmp(u_int *, u_int *);


u_int           prt_tol, p;


/* ************************************************************************** */


#define show_perm(PP) ((PP->size < prt_tol) ? out_perm(PP) :\
                      printf("Permutation: size: %d\n", PP->size))
#define show_vec(VV) ((VV->dim < prt_tol) ? out_vec(VV) :\
                     printf("Vector: dim: %d\n", VV->dim))
#define show_mat(MM) ((MM->m < prt_tol) ? out_mat(MM) :\
                     printf("Matrix: m: %d by 2\n", MM->m))






PERM *decomp(sp_mat *L, PERM *P, PERM *A, PERM *B, PERM *S, int rec_lvl)
{
        int             i, j, k, l, n = L->m, tempA, tempB, tempH,
                        inA, inB, inA1, inB1, inD, inS, inAdash, inBdash,
                        N_nonzero_v, max_degree, vmax, imax, ev_good_enough;
        double          medval, L2old, L2new, L2tol = 1e-04, Rtol = 1e-01,
                        yi, yj, sum, alpha, beta; 
        PERM            *Adash, *Bdash, *A1, *B1, *AA, *AB, *AS, *BA, *BB, *BS,
                        *pivot;
        VEC             *a, *b, *c, *d, *x0, *w, *y, *z, *resid, *V, *W, *TMP;
        MAT             *C, *I, *Q, *T, *D, *H;
        row_elt         *elt_list;
        sp_mat          *AL, *BL;


        if (!L || !A || !B || !S || !P)
                error(E_NULL, "decomp");
        if (L->m != L->n || P->size != n)
                error(E_SIZES, "decomp");

        printf("At top of decomp. rec_lvl: %d\n", rec_lvl);

/* 0. Set up all the required matrix and vector elements. Firstly, initialise
all those with dimensions fixed, depending on "n".*/

        Adash = get_perm(n);                            Bdash = get_perm(n);
        A1 = get_perm(n);                               B1 = get_perm(n);
        y = get_vec(n);                                 x0 = get_vec(n);
        resid = get_vec(n);                             TMP = get_vec(n);
        V = get_vec(n);                                 W = get_vec(n);
        D = get_mat(n*p, 2);                            H = get_mat(n*p, 2);

/* Initialise size of variable-size data-structures 1. They are soon resized. */

        a = get_vec(1);                                 b = get_vec(1);
        c = get_vec(1);                                 d = get_vec(1);
        I = get_mat(1, 1);                              T = get_mat(1, 1);
        C = get_mat(1, 1);                              Q = get_mat(1, 1);
        w = get_vec(1);                                 z = get_vec(1);        
        pivot = get_perm(1);


/* 1. Use the Lanczos method and a tridiagonal eigendecomposition routine to
calculate the 2nd eigenvalue and corresponding eigenvector of L. Pothen, et al.
suggest the choice for the initial x0, and we begin with W = normalised(x0) and
V = LW. */
        for (i = 0; i < n; x0->ve[i] = i - (n - 1)/2, i++);
        L2old = 0;                                      ev_good_enough = 0;
        while (!ev_good_enough)
        {
                j = 0;
                sv_mlt(1.0/n2(x0), x0, W);
                sp_mv_mlt(L, W, V);
                L2new = L2old + 2*L2tol;                beta = 1;
                while (j < 2 || fabs((L2old - L2new)/L2new) > L2tol)
                {
                        j++;
                        if (j*n > 250000)
                        {
                                printf("j.n = %d.%d = %d ", j, n, j*n);
                                printf("Using too much memory, cutting out!\n");
                                exit(0);
                        }
                        a = v_resize(a, j);             b = v_resize(b, j-1);
                        if (j > 1)
                                b->ve[j-2] = beta;
                        Q = m_resize(Q, n, j);          set_col(Q, j-1, W);
/* Store W in Q. */
                        a->ve[j-1] = alpha = in_prod(W, V);
                        v_mltadd(V, W, -alpha, V);


/* Orthogonalise V relative to e. */
                        for (i = sum = 0; i < n; sum += V->ve[i], i++);
                        for (i = 0, sum = sum/n; i < n; V->ve[i] -= sum, i++);

                        beta = n2(V);                   cp_vec(W, TMP);
                        sv_mlt(1/beta, V, W);           sv_mlt(-beta, TMP, V);
                        sp_mv_mlt(L, W, TMP);           v_add(V, TMP, V);

                        c = v_resize(c, j);             d = v_resize(d, j-1);
                        c = cp_vec(a, c);               d = cp_vec(b, d);
/* trieig(c, d, M) takes a tridiagonal matrix with diagonal entries c and
off-diagonal entries d, finds the eigenvalues, and stores then in c. */
                        trieig(c, d, MNULL);
                        L2old = L2new;                  L2new = select(c, 1);
                }
                if (L2new < 0)
                {
                        printf("Negative ev!\n");       exit(0);
                }

/* Set up the T and I matrices, and the w vector. */
                T = m_resize(T, j, j);                  zero_mat(T);
                for (i = 0; i < j - 1; i++)
                {
                        T->me[i][i] = a->ve[i];
                        T->me[i + 1][i] = T->me[i][i + 1] = b->ve[i];
                }
                T->me[j - 1][j - 1] = a->ve[j - 1];
                I = m_resize(I, j, j);                  I = id_mat(I);
                pivot = px_resize(pivot, j);            pivot = px_id(pivot);
                w = v_resize(w, j);                     rand_vec(w);
                z = v_resize(z, j);                     C = m_resize(C, j, j);

                printf("\nFinding y: j = %d\n\tlambda2\t\tn2(resid)\n", j);

                C = sm_mlt(-L2new, I, C);               C = m_add(T, C, C);
                for (i = 0; i < C->m; i++)
                        if (C->me[i][i] == 0)
                                C->me[i][i] = MACHEPS;
                C = LUfactor(C, pivot);
                z = LUsolve(C, pivot, w, z);
                w = sv_mlt(1.0/n2(z), z, w);
                L2old = L2new;
                L2new = in_prod(w, mv_mlt(T, w, z));

/* If the residual of eigenvector/value pair || L.y - L2.y ||_2 / || y ||_2 is
too large, run again ... */
                y = mv_mlt(Q, w, y);
                resid = sp_mv_mlt(L, y, resid);
                resid = v_mltadd(resid, y, -L2new, resid);
                printf("\t%25.20g\t\t%25.20g\n", L2new, n2(resid));
                ev_good_enough = (n2(resid) < Rtol*n2(y));
                if (!ev_good_enough)
                        x0 = cp_vec(y, x0);
        }

        printf("Using lambda2 = %20.15g\n", L2new); 
        printf("y\, the 2nd eigenvector of L: "); show_vec(y);
        
/* Most of the continuous mathematics data structures are no longer needed. */
        freevec(a);     freevec(b);     freevec(c);     freevec(d);
        freevec(x0);    freevec(w);     freevec(z);     freeperm(pivot);
        freevec(V);     freevec(W);     freevec(TMP);
        freemat(C);     freemat(I);     freemat(Q);     freemat(T); 








/* 2. Calculate Adash, Bdash and H. Adash and Bdash are the vertex sets
generated by the partitioning of the vertices of G by H. H is an edge separator
set of G, found using Fiedler's method. y is the 2nd eigenvector of the
Laplacian matrix of G, with median value medval. */

/* 2.1. Calculate the median of the elements of the 2nd eigenvector. */
        medval = select(y, (u_int) ((n + 1)/2));

/* 2.2.1.  Set up Adash and Bdash. */

        for (i = Adash->size = Bdash->size = 0; i < y->dim; i++)
                if (y->ve[i] <= medval)
                        Adash->pe[Adash->size++] = i;
                else
                        Bdash->pe[Bdash->size++] = i;

        printf("First Pass Adash and Bdash made\n");
        show_perm(Adash);                               show_perm(Bdash);
        
/* 2.2. Set up H. Search through the upper half triangle of L, and, for
each edge encountered, insert it into H only if the eigenvaluation of
the vertices crosses the median. */

        for (i = H->m = 0; i < L->m; i++)
        {
                elt_list = (L->row[i]).elt;
                yi = y->ve[i];
                for (k = 0; k < (L->row[i]).len; k++)
                        if ((j = elt_list[k].col) < i)
                        {
                                yj = y->ve[j];
                                if (yi <= medval && yj > medval)
                                {
                                        H->me[H->m][0] = i;
                                        H->me[H->m++][1] = j;
                                }
                                else if (yj <= medval && yi > medval)
                                {
                                        H->me[H->m][0] = j;
                                        H->me[H->m++][1] = i;
                                }
                        }
        }
         printf("First Pass H made\n");                 show_mat(H);

/* 2.2.2. If |Adash| - |Bdash| > 1, move enough vertices with components equal
to medval from Adash to Bdash. Not easy! Also must correct H. */

        for (i = 0; i < Adash->size && Adash->size - Bdash->size > 1; i++)
                if (y->ve[i] == medval)
                {
                        for (j = 0; Adash->pe[j] != i; j++);
                        Adash->pe[j] = Adash->pe[--Adash->size];
                        Bdash->pe[Bdash->size++] = i;
                        for (k = 0; k < H->m; k++)
                                for (l = 0; H->me[k][0] == i && l < 2; l++)
                                        H->me[k--][l] = H->me[H->m--][l];
                }
        freevec(y);

/* Sort Adash and Bdash. */
        qsort(Adash->pe, Adash->size, sizeof(u_int), icmp);
        qsort(Bdash->pe, Bdash->size, sizeof(u_int), icmp);

        printf("2nd Pass Adash: ");                     show_perm(Adash);
        printf("\n2nd Pass Bdash: ");                   show_perm(Bdash);
        printf("\n2nd Pass H: ");                       show_mat(H);





/* 3. Calculate S, a vertex separator set of G. A reminder that H is an edge
separator set, and Adash and Bdash are a disjoint cover of the vertex set of G.
A1 and B1 are the respective subsets of Adash and Bdash that consist of only
vertices that adjoin edges in H. The two vertex sets A and B are found, such
that A = Adash \ S, B = Bdash \ S). Respective edges are not relevant to the
further progress of "decomp", and so are not found. S is optimally as small as
possible, and one is found that is hopefully relatively small, by the better of
2 algorithms that seem like a quick and good way:
        1. Collect vertices of reducing degree from A1 U B1 until H is covered.
        2. Use all the vertices of the smaller of the two sets, A1 and B1. */

/* 3.1.1. Set up D, a listing of the degrees of vertices in H. */

        for (i = D->m = 0; i < H->m; i++)
                for (j = 0; j < 2; j++)
                {
                        tempH = (u_int) H->me[i][j];
                        for (k = inD = 0; k < D->m; k++)
                                if ((u_int) D->me[k][0] == tempH)
                                        D->me[k][inD = 1]++;
                        if (!inD)
                        {
                                D->me[D->m][0] = tempH;
                                D->me[D->m++][1]++;
                        }
                }

        printf("\nD: ");                                show_mat(D);
                        
/* 3.1.2. Set up A1 and B1, respective sides of H. */
        
        for (i = A1->size = B1->size = 0; i < H->m; i++)
        {
                tempH = H->me[i][0];
                for (k = inA1 = 0; k < A1->size; k++)
                        inA1 = inA1 || A1->pe[k] == tempH;
                if (!inA1)
                        A1->pe[A1->size++] = tempH;
                tempH = H->me[i][1];
                for (k = inB1 = 0; k < B1->size; k++)
                        inB1 = inB1 || B1->pe[k] == tempH;
                if (!inB1)
                        B1->pe[B1->size++] = tempH;
        }
/* Sort A1 and B1. Not currently useful, but good for the future. */
        qsort(A1->pe, A1->size, sizeof(u_int), icmp);
        qsort(B1->pe, B1->size, sizeof(u_int), icmp);

        printf("2nd Pass A1: ");                        show_perm(A1);
        printf("\n2nd Pass B1: ");                      show_perm(B1);

/* 3.2. Establish the permutation of the vertex separator set S, via a grungy
but apparently working piece of code. The algorithm used:

        1. Find the vertex of highest degree in D, such that:
                D->me[imax][0] = vmax; D->me[imax][1] = max_degree;
        2. Set the degree of this vertex to zero, decrement the number of
        non-zero vertices, and add this vertex to S.
                D->me[imax][1] = 0; N_non_zero_v--; S->pe[S->size++] = vmax;
        3. Search out the edges in H containing this vertex, and decrement the
        vertices adjacent to it in D.
        4. Repeat until no non-zero degree vertices remain. */

        N_nonzero_v = D->m;                             S->size = 0;
        while (N_nonzero_v > 0)
        {
                for (i = max_degree = 0; i < D->m; i++)
                        if (D->me[i][1] > max_degree)                
                                max_degree = (u_int) D->me[imax = i][1];

                D->me[imax][1] = 0;                     N_nonzero_v--;
                S->pe[S->size++] = vmax = (u_int) D->me[imax][0];

                for (i = 0; i < H->m && N_nonzero_v > 0 && max_degree > 0; i++)
                        for (j = 0; j < 2 && N_nonzero_v > 0 && max_degree > 0; j++)
                                if ((u_int) H->me[i][j] == vmax)
                                {
                                        k = (u_int) H->me[i][!j];
                                        for (l = 0; l < D->m && (u_int) D->me[l][0] != k; l++);
                                        if (D->me[l][1] > 0)
                                        {
                                                D->me[l][1]--;
                                                max_degree--;
                                                if ((u_int) D->me[l][1] == 0)
                                                        N_nonzero_v--;
                                        }
                                }
        }
/* Sort S. */
        qsort(S->pe, S->size, sizeof(u_int), icmp);

        printf("\nFirst pass S: ");                     show_perm(S);

/* 3.3. Test to see if the S that has been created is larger than the
smaller of A1 and B1. If it is, replace S with this smaller set. */
        
        if (S->size > A1->size)
        {
                S = cp_perm(A1, S);                     S->size = A1->size;
        }
        if (S->size > B1->size)
        {
                S = cp_perm(B1, S);                     S->size = B1->size;
        }
        printf("\nSecond pass S: ");                    show_perm(S);

/* 3.4. Create A, B, such that A, B and S are a disjoint cover of N. */
        for (i = A->size = 0; i < Adash->size; i++)
        {
                k = Adash->pe[i];
                for (j = inS = 0; j < S->size; j++)
                        inS = inS || S->pe[j] == k;
                if (!inS)
                        A->pe[A->size++] = k;
        }
        for (i = B->size = 0; i < Bdash->size; i++)
        {
                k = Bdash->pe[i];
                for (j = inS = 0; j < S->size; j++)
                        inS = inS || S->pe[j] == k;
                if (!inS)
                        B->pe[B->size++] = k;
        }
/* Sort A and B. */
        qsort(A->pe, A->size, sizeof(u_int), icmp);
        qsort(B->pe, B->size, sizeof(u_int), icmp);


/* Print out the vital statistics of the whole process. */
        printf("\nSizes of various elements:\n\n");
        printf("\t# Adash, Bdash:   %d\t%d\n", Adash->size, Bdash->size);
        printf("\tLambda2:          %f\n", L2new);
        printf("\t# D:              %d\n", D->m);
        printf("\t# H:              %d\n", H->m);
        printf("\t# A1 and B1:      %d\t%d\n", A1->size, B1->size);
        printf("A: ");                                  show_perm(A);
        printf("B: ");                                  show_perm(B);
        printf("S: ");                                  show_perm(S);
        printf("*******************************************\n\n");

/* Dump unneeded data-structures here! */
        freeperm(A1);           freeperm(Adash);        freemat(D);
        freeperm(B1);           freeperm(Bdash);        freemat(H);


/* 4. Recursively partition the sets A and B by calling "decomp" on them. */

/* 4.1. Decompose A and B if their sizes are > 3. For example, for A, a sparse
matrix AL is created (found by extracting the parts of L with indices in A),
then decomp is called with AL as L, A as P, and dummies AA, AB, AS as the other
parameters. In fact, the routine does not actually need these other parameters,
but they are used, as they may be required for further programming to also
return the structure of the final permutation. Required structures are created
in situ. */
        if (A->size > 3 && rec_lvl != -1)
        {
                AA = get_perm(A->size);
                AB = get_perm(A->size);
                AS = get_perm(A->size);
                AL = sp_get_mat(A->size, A->size, 3*p);

/* The following 9 lines of _slow_ code (to set the entries of AL) are a waste
of time, but I haven't yet figured out what else to do. Ditto for BL. In fact,
it appears that this is the slowest part of the whole algorithm! */

                for (i = 0; i < A->size; i++)
                        for (j = 0; j < A->size; j++)
                        {
                                k = (u_int) sp_get_val(L, A->pe[i], A->pe[j]);
                                if (k != 0)
                                        sp_set_val(AL, i, j, k);
                        }
                for (i = 0; i < A->size; i++)
                        sp_set_val(AL, i, i, (AL->row[i]).len - 1);

/* *****************************************************************************
The alternative version uses the elements of A as sorted into increasing order.
Now we can scan through the rows of L that corresponded to entries in A, and
know that we only have to check the entries in the row of L up until the column
number becomes >= row number. It looks something like:

                for (i = 0; i < A->size; i++)
                {
                        k = A->pe[i];
                        elt_list = (L->row[k]).elt;
                        for (j = 0; j < (L->row[k]).len; j++)
                        {
                                for (l = inA = 0; l < A->size; l++)
                                        inA = inA || A->pe[l] == (elt_list[j]).val;
                                if (inA)
                                        sp_set_val(AL, i, j, -1);
                        }
                }
                for (i = 0; i < A->size; i++)
                        sp_set_val(AL, i, i, (AL->row[i]).len);
***************************************************************************** */


                if (!ck_symm(AL))
                {
                        printf("Quitting as AL not symmetric!\n");
                        exit(0);
                }
                if (!ck_sums(AL))
                {
                        printf("Quitting as AL has invalid row sums!\n");
                        exit(0);
                }
                decomp(AL, A, AA, AB, AS, ++rec_lvl);           rec_lvl--;
                freeperm(AA);   freeperm(AB);   freeperm(AS);   sp_free_mat(AL);
        }





        if (B->size > 3 && rec_lvl != -1)
        {
                BA = get_perm(B->size);
                BB = get_perm(B->size);
                BS = get_perm(B->size);
                BL = sp_get_mat(B->size, B->size, 3*p);

                for (i = 0; i < B->size; i++)
                        for (j = 0; j < B->size; j++)
                        {
                                k = (u_int) sp_get_val(L, B->pe[i], B->pe[j]);
                                if (k != 0)
                                        sp_set_val(BL, i, j, k);
                        }
                for (i = 0; i < B->size; i++)
                        sp_set_val(BL, i, i, (BL->row[i]).len - 1);

                if (!ck_symm(BL))
                {
                        printf("Quitting as BL not symmetric!\n");
                        exit(0);
                }
                if (!ck_sums(BL))
                {
                        printf("Quitting as BL has invalid row sums!\n");
                        exit(0);
                }
                decomp(BL, B, BA, BB, BS, ++rec_lvl);           rec_lvl--;
                freeperm(BA);   freeperm(BB);   freeperm(BS);   sp_free_mat(BL);
        }

/* 4.2. Map the apparent A, B, S to the actual A, B, S, using P, the permutation
of their actual names. */
        for (i = 0; i < A->size; i++)
                A->pe[i] = P->pe[A->pe[i]];
        for (i = 0; i < B->size; i++)
                B->pe[i] = P->pe[B->pe[i]];
        for (i = 0; i < S->size; i++)
                S->pe[i] = P->pe[S->pe[i]];

        for (i = 0; i < A->size; i++)
                P->pe[i] = A->pe[i];
        j = i;
        for (i = 0; i < B->size; i++)
                P->pe[i + j] = B->pe[i];
        j += i;
        for (i = 0; i < S->size; i++)
                P->pe[i + j] = S->pe[i];

        return P;
}


/* ************************************************************************** */


int icmp(u_int *p1, u_int *p2)
{
        return *p1 - *p2;
}
\end{verbatim} 
\normalsize\rm

\pagebreak

\subsection{Listing of \texttt{mk\_sp\_graph}}
\scriptsize\tt 
\begin{verbatim}
/*                            mk_sp_graph                             */


/* *********************************************************************


Contents  : This file is called "mk_sp_graph.c", and contains the
	    function "mk_sp_graph".

     Aim  : Generate a random G, the sparse matrix of a (sparse)
	    undirected, unvaluated graph on n vertices with an average
	    degree of p (number of edges/vertex). Further, the graph is
	    almost planar, and this is achieved by using only a narrow
	    bandwidth (q) in the original matrix, followed by a
	    scrambling of its indices.

 Language : ANSI Standard C

  Author  : David  De Wit    March 5  -  May 27  1991


********************************************************************* */


#include <stdio.h>
#include "matrix.h"
#include "sparse.h"


sp_mat *mk_sp_graph(sp_mat *G, u_int p, u_int q)
{
        u_int           i, j, n = G->m;
        double          temp, limit;

        
        srand(4);
        limit = ((double) p)*2147483648/((double) n);
        for (i = 0; i < n; i++)
                for (j = i + 1; (j < n) && (j < i + q); j++)
                        if (rand() < limit)
                        {
                                 temp = sp_set_val(G, i, j, 1.0);
                                temp = sp_set_val(G, j, i, 1.0);
                        }
                
        return G;
}
\end{verbatim} 
\normalsize\rm

\pagebreak

\subsection{Listing of \texttt{select}}
\scriptsize\tt 
\begin{verbatim}
/*                            select                                  */


/* *********************************************************************


Contents  : This file is called "select.c", and contains the function
"select".

     Aim  : Select the kth smallest entry in vector a. On exit, the
	    desired element is in its correct place. In particular,
	    calling select(a, int((a->dim + 1)/2)) finds the median.

            See "Algorithms", Sedgewick (1983), p128. QA76.6.S435

 Language : ANSI Standard C

  Author  : David  De Wit    March 4  -  May 27  1991


********************************************************************* */


#include <stdio.h>
#include "matrix.h"


double select(VEC *a, u_int k)
{
        int             i, j, l, r;
        double          t, v;
        VEC             *b;


        b = get_vec(a->dim);            b = cp_vec(a, b);
        l = 0;                          r = b->dim - 1;
        while (r > l)
        {
                v = b->ve[r];           i = l - 1;      j = r;
                do
                {
                        for (i++; b->ve[i] < v; i++);
                        for (j--; b->ve[j] > v; j--);
                        t = b->ve[i];
                        b->ve[i] = b->ve[j];
                        b->ve[j] = t;
                } while (j > i);
                b->ve[j] = b->ve[i];
                b->ve[i] = b->ve[r];
                b->ve[r] = t;
                if (i >= k)
                        r = i - 1;
                if (i <= k)
                        l = i + 1;
        }
        return (b->ve[k - 1]);
}
\end{verbatim} 
\normalsize\rm

\pagebreak

\subsection{Listing of \texttt{gr\_lap}}
\scriptsize\tt 
\begin{verbatim}
/*                            gr_lap                                  */
 

/* *********************************************************************


Contents  : This file is called "gr_lap.c", and contains the program of
	    the same name.

     Aim  : Make the graph for solving Laplace's Equation on an M x N
	     grid.

 Language : ANSI Standard C

  Author  : David  De Wit (and David Stewart)  May 8  -  May 27  1991


********************************************************************* */


#include <stdio.h>
#include <strings.h>
#include "matrix.h"
#include "sparse.h"


#define index(i,j) (N*((i)-1)+(j)-1)


main(int argc, char *argv[])
{
        u_int           i, j, M = atoi(argv[1]), N = atoi(argv[2]);
        char            *outname;
        sp_mat          *A;


        A = sp_get_mat(M*N, M*N, 5);
        for (i = 1; i <= M; i++)
                for (j = 1; j <= N; j++)
                {
                        if (i < M)
                                sp_set_val(A, index(i,j), index(i+1,j), 1);
                        if (i > 1)
                                sp_set_val(A, index(i,j), index(i-1,j), 1);
                        if (j < N)
                                sp_set_val(A, index(i,j), index(i,j+1), 1);
                        if (j > 1)
                                sp_set_val(A, index(i,j), index(i,j-1), 1);
                }
        outname = strcat(strcat(strcat("Lap.", argv[1]), "."), argv[2]);
        sp_fout_mat(fopen(outname, "w"), A);
}
\end{verbatim} 
\normalsize\rm

\pagebreak

\subsection{Listing of \texttt{testdc}}
\scriptsize\tt 
\begin{verbatim}
/*                            testdc                                  */
 

/* *********************************************************************


Contents  : This file is called "testdc.c", and contains the program of
	    the same name. It calls the function "decomp".

     Aim  : Test the function "decomp" by setting-up and solving a
	     problem.

 Language : ANSI Standard C

  Author  : David  De Wit    March 5  -  June 19  1991


********************************************************************* */


#include <stdio.h>
#include "matrix.h"
#include "matrix2.h"
#include "sparse.h"


extern u_int    prt_tol, p, q;
extern PERM     *decomp(sp_mat *, PERM *, PERM *, PERM *, PERM *, int);
extern sp_mat   *mk_sp_graph(sp_mat *, u_int, u_int);
extern int      ck_symm(sp_mat *), ck_sums(sp_mat *);


main(int argc, char *argv[])
{
        u_int           i, j, k, n = atoi(argv[1]), q = 6, nfiles = 11;
        int             j_idx, rec_lvl = atoi(argv[4]);
        PERM            *A, *B, *S, *P;
        FILE            *fp;
        row_elt         *e;
        sp_row          *r;
        sp_mat          *G;
        char            *fname[] = {"Idit", "Moshe", "Itzchak", "Shimuel",
                                    "Ishmail", "Yacov", "Yair", "Arieh",
                                    "Aaron", "Schlomo", "Shimshon"};


/*. 0. Initialise some structures and constants for the problem. */
        prt_tol = atoi(argv[3]);                p = atoi(argv[2]);
        A = get_perm(n);                        B = get_perm(n);
        S = get_perm(n);                        P = get_perm(n);
        P = px_id(P);
                


















/* 1. Set up problem by randomly generating or reading in a sparse matrix. */
        printf("\nEnter a choice for the initial graph:\n\n");
        printf("\t0: Make a random sparse graph\n");
        for (i = 0; i < nfiles; i++)
                printf("\t%1d: Read file \"%s\"\n", i + 1, fname[i]);
        printf("\nYour Choice:\n");             scanf("%d", &i);
        if (i)
        {        
                printf("\nReading a sparse graph from \"%s\"\n", fname[i - 1]);
                fp = fopen(fname[i - 1], "r");
                G = sp_fin_mat(fp);
        }
        else
        {
                printf("\nMaking a sparse graph on %d vertices\, ", n);
                printf("at an average %d edges\/vertex.\n", p);
                G = sp_get_mat(n, n, 3*p);
                G = mk_sp_graph(G, p, q);
        }        

        if (G->m != n)
                error(E_SIZES, "testdc");


/* 2. Calculate the sparse matrix of the graph Laplacian (L) of the graph
represented by the sparse matrix G. Defining d[i] as the degree of vertex i in
G; then L = diag(d) - G. All row and column sums in L = 0. G overwrites L. */

/* Ensure diagonal entries are in rows. */
        for (i = 0; i < G->m; sp_set_val(G, i, i, 1.0), i++);

/* Set the values of the entries. */
        for (i = 0; i < G->m; i++)
        {
                r = &(G->row[i]);
                /* scan entries of row r */
                for (j_idx = 0, e = r->elt; j_idx < r->len; j_idx++, e++)
                    if (e->col == i)
                        e->val = r->len - 1;        /* diagonal entry */
                    else
                        e->val = -1.0;              /* off-diagonal entry */
        }

        if (!ck_symm(G))
        {
                printf("Quitting as graph Laplacian not symmetric!\n");
                exit(0);
        }
        if (!ck_sums(G))
        {
                printf("Quitting as graph Laplacian has invalid row sums!\n");
                exit(0);
        }


/* 3. Call the function to do the partitioning. */
        printf("\nCalling decomp from testdc ...\n\n");
        P = decomp(G, P, A, B, S, rec_lvl);
        printf("\nBack in Kansas ... P:\n");
        if (P->size < prt_tol)
                out_perm(P);
        else
                printf("Permutation, size: %d\n", P->size);
}








/* ************************************************************************** */


int ck_symm(sp_mat *L)
{
        int             i;
        static VEC      *x=VNULL, *y1=VNULL, *y2=VNULL;


        if (!L)
                error(E_NULL,"ck_symm");
        if (L->m != L->n)
                return FALSE;

        x  = v_resize(x, L->m);                 y1 = v_resize(y1, L->m);
        y2 = v_resize(y2, L->m);

        for (i = 0; i < 3; i++)
        {
                rand_vec(x);                    y1 = sp_mv_mlt(L, x, y1);
                y2 = sp_vm_mlt(L, x, y2);       y1 = v_sub(y1, y2, y1);
                if (n2(y1) > L->m*MACHEPS)
                        return FALSE;
        }

        return TRUE;
}


/* ************************************************************************** */


int ck_sums(sp_mat *L)
{
        int             i, j;
        double          sum;
        sp_row          *r;


        if (!L)
                error(E_NULL,"ck_sums");
        for (i = 0; i < L->m; i++)
        {
                r = &(L->row[i]);
                for (j = sum = 0; j < r->len; j++)
                        sum += r->elt[j].val;
                if (fabs(sum) > L->m*MACHEPS)
                        return FALSE;
        }

        return TRUE;
}
\end{verbatim} 
\normalsize\rm

\pagebreak

%%%%%%%%%%%%%%%%%%%%%%%%%%%%%%%%%%%%%%%%%%%%%%%%%%%%%%%%%%%%%%%%%%%%%%%%

\bibliographystyle{plain}
\bibliography{DeWit91}

\end{document}